\documentclass[review]{elsarticle}

\usepackage{lineno,hyperref}
\usepackage{amssymb,amsfonts,amsmath,euscript}
\usepackage{amsmath, amssymb,graphics,amsfonts, array, bm}
\numberwithin{equation}{section}
\modulolinenumbers[5]

\begin{document}

\begin{frontmatter}
\title{Sufficient conditions for convergence of multiple Fourier series with $J_k$-lacunary
sequence of rectangular partial sums in terms of Weyl multipliers}
\author[]{I.L. Bloshanskii\corref{mycorrespondingauthor}}
\author[]{S.K. Bloshanskaya}
\author[]{D.A. Grafov}

\begin{abstract}
We obtain sufficient conditions for convergence (almost everywhere) of multiple trigonometric Fourier series of functions $f$ in $L_2$
in terms of Weyl multipliers.
We consider the case where rectangular partial sums
of Fourier series $S_n(x;f)$ have indices
$n=(n_1,\dots,n_N) \in \mathbb Z^N$, $N\ge 3$, in which $k$ $(1\leq k\leq N-2)$ components on the places
$\{j_1,\dots,j_k\}=J_k \subset \{1,\dots,N\} = M$
are elements of (single) lacunary sequences
(i.e., we consider the, so called, multiple Fourier series with $J_k$-lacunary sequence of partial sums).
We prove that for any sample $J_k\subset M$ the Weyl multiplier for convergence of these series
has the form $W(\nu)=\prod \limits_{j=1}^{N-k} \log(|\nu_{{\alpha}_j}|+2)$, where $\alpha_j\in M\setminus J_k $, $\nu=(\nu_1,\dots,\nu_N)\in{\mathbb Z}^N$.
So, the "one-dimensional"\, Weyl multiplier -- $\log(|\cdot|+2)$--
presents in $W(\nu)$ only on the places of "free"\, (nonlacunary) components of the vector $\nu$.
Earlier, in the case where $N-1$ components of the index $n$ are elements of lacunary sequences, convergence almost everywhere for multiple Fourier series was obtained in 1977 by M.\,Kojima in the classes $L_p$, $p>1$, and by D.\,K.\,Sanadze, Sh.\,V.\,Kheladze in Orlizc class. Note, that presence of two or more "free"\, components in the index $n$ (as follows from the results by Ch.\,Fefferman (1971)) does not guarantee the convergence almost everywhere of $S_n(x;f)$ for $N\geq 3$ even in the class of continuous functions.
\end{abstract}
\end{frontmatter}

\section{Introduction}\label{s2}

{\bf 1.} Consider the $N$-dimensional Euclidean space
$\mathbb R^N$, whose elements
will be denoted as $x=(x_1,\dots,x_N)$, and set
$(nx)=n_1x_1+\dots+n_Nx_N$.
We introduce ${\mathbb R}^N_{\sigma}=\{(x_1,\,\dots\,,x_N)  \in {\mathbb R}^N: x_j\geq
\sigma,\ j=1,\dots,N\}$, $\sigma\in {\mathbb R}^1$, and the set ${\mathbb
Z}^N\subset{\mathbb R}^N$ of all vectors
with integer coordinates. Set ${\mathbb Z}^N_{\sigma} = {\mathbb R}^N_{\sigma} \cap
{\mathbb Z}^N$.

Let a $2\pi$-periodic (in each argument) function $f\in
L_1({\mathbb T}^N)$, where ${\mathbb T}^N=\{x\in{{\mathbb R}^N}: -\pi \leq x_j<\pi, j=1,\dots,N\}$, be expanded in a
multiple trigonometric Fourier series: $f(x)\thicksim\sum_{\nu\in {\mathbb Z}^N}c_\nu e^{i(\nu x)}.$

For any vector $n=(n_1,\dots,n_N)\in{\mathbb Z}_0^N $ consider a rectangular partial sum of these series
\begin{equation}\label{0.1}
S_n(x;f)=\sum_{|\nu_1|\leq n_1}\dots\sum_{|\nu_N|\leq n_N}c_\nu e^{i(\nu x)}.
\end {equation}

The main purpose of our investigation is to study the behavior on $\mathbb T^N$ of the partial sum \eqref{0.1} as  $n\to \infty$ (i.e. $\min\limits_{1 \leq j \leq N}n_j\to\infty$), depending on the restrictions imposed as on the function $f$, so as on the components $n_1,\dots,n_N$ of the vector $n$ -- the index of $S_n(x;f)$.

In 1971  P.\,Sjolin \cite{Sjolin} proved that for any lacunary sequence\,\footnote{\,A sequence $\{n^{(s)}\}$,
$n^{(s)}\in {\mathbb Z}^1_1$, is called lacunary, if $n^{(1)}=1$ and $\frac
{n^{(s+1)}}{n^{(s)}}\geq q>1$, $s=1,2,\dots \enskip.$}
$\{n_1^{(\lambda_1)}\}, \ n_1^{(\lambda_1)}\in {\mathbb Z}_1^1,
\lambda_1=1,2,\dots,$ and for any function $f\in L_p({\mathbb T}^2)$, $p>1$,
$$
\lim\limits_{\lambda_1, \, n_2\to \infty} S_{n_1^{(\lambda_1)},\,n_2}(x;f)=
f(x) \quad\text {almost everywhere (a.e.) on}\quad {\mathbb T}^2.\ \footnote{\,In 1970 N. R. Tevzadze \cite{Tevzadze} obtained the following result: for any given two sequences of integers $\{n_j^{(l)}\}$, $j=1,2$, increasing to $\infty$, $n_j^{(l)}\in {\mathbb Z}_1^1$, $\l=1,2,\dots$, $S_{n_1^{(l)},\,n_2^{(l)}}(x;f)$ converges to $f(x)$ a.e. on ${\mathbb T}^2$ for $f\in L_2({\mathbb T}^2)$.}
$$

In 1977 M.\,Kojima \cite{Kojima} generalized P.\,Sjolin's result by
proving that, if a function  $f \in L_p({\mathbb T}^N)$, $p> 1$, $N
\geq 2$, and $\{n_j^{(\lambda_j)}\}$, $n_j^{(\lambda_j)}\in {\mathbb Z}_1^1,
\lambda_j=1,2,\dots, j=1,\dots,N-1$, are lacunary sequences, then
$$
\lim\limits_{\lambda_1,\, \dots,\, \lambda_{N-1},\, n_N \to \infty}
S_{n_1^{(\lambda_1)}, \dots,\, n_{N-1}^{(\lambda_{N-1})},\,n_N}(x;f) = f(x)
\quad\text {a.e. on}\quad {\mathbb T}^N.\footnote{\,In the same 1977 for $f \in L(\log^+ L)^{3N-2}(\mathbb T^N)$ analogous result was obtained by D.\,K.\,Sanadze, Sh.\,V.\,Kheladze \cite{Sanadze}.}
$$

However, as soon as we remain "free"\, two components of the vector $n=(n_1,\dots, n_N)
\in \mathbb Z^N_0$ -- the index of $S_n(x;f)$ (in particular, in the case where they are not elements of any lacunary sequences), even the class of continuous functions ${\mathbb C}(\mathbb T^N)$, $N\geq3$, does not remain the "class of convergence a.e."\, of the considered expansions; this can be easily shown using Ch.\,Fefferman's function from \cite{Fefferman} (see, e.g. \cite[Theorem 2]{Kojima}). Nevertheless, some conditions can be imposed on the "nonlacunary"\, components of the vector $n$ (in the sequence of indices of partial sums), such that even the classes $L_p({\mathbb T}^N), p>1$, for $N\geq3$, remain the classes of convergence a.e. (of the considered expansions), in the case where there are more than one nonlacunary components; moreover, for the certain subsets of $L_p({\mathbb T}^N)$ all nonlacunary components can be even "free".

For functions in the classes $L_p({\mathbb T}^N), p>1, N\ge 3$, I.\,L.~Bloshanskii and D.\,A.~Grafov \cite{Bloshanskii1} proved convergence a.e. on ${\mathbb T}^N$ of the sequence of partial sums of multiple trigonometric Fourier series whose indices $n$ contain $k$ lacunary components, $1\leq k\leq N-2$, while the rest $N-k$ nonlacunary components $n_j$ of the vector $n$ satisfy restrictions: $c_1\le n_j/n_m\le c_2,$  where $c_1, c_2= const$ (so, along the nonlacunary components, the summation over an extending system of rectangles takes place).

A question naturally arises to find such classes of functions which guarantee convergence a.e. of the sequence of partial sums of multiple Fourier series with indices $n$ whose $k$ components are lacunary ($1\leq k\leq N-2$), and at the same time, $N-k$ nonlacunary components of the vector $n$ are either "more free"\, than in the paper \cite{Bloshanskii1}, or "free at all"\,.

{\bf 2.} In the present paper we give an answer to this question in terms of Weyl multipliers.

{\bf Definition.} {\it A sequence $W(\nu)$, $\nu\in{\mathbb Z}^N$, is called a Weyl multiplier for rectangular convergence of a multiple trigonometric Fourier series if, first, it satisfies the conditions:

$1.$ $W(\nu)>0$,\ $\nu\in{\mathbb Z}^N$;

$2.$ $W(\nu_1,\dots,\nu_N)=W(|\nu_1|,\dots,|\nu_N|)$,\quad $\nu\in{\mathbb Z}^N$;

$3.$ $W(\nu_1,\dots, \nu_{j-1},\nu_j+1,\nu_{j+1},\dots, \nu_N)\geq W(\nu_1,\dots, \nu_{j-1},\nu_j,\nu_{j+1},\dots,\nu_N)$,\quad $j=1,\dots,N$,\quad $\nu\in{\mathbb Z}_0^N$;

and, second, if the convergence of the series $\sum_{\nu\in {\mathbb Z}^N}|c_\nu(f)|^2W(\nu)<+\infty$ implies that the Fourier series of the function $f\in L_2(\mathbb T^N)$ converge over rectangles a.e. on $\mathbb T^N$}.

From the L.\,Carleson theorem \cite{Carleson} it follows that in the one-dimensional case the Weyl multiplier $W(\nu)=1$, $\nu\in{\mathbb Z}^1$. For the $N$-multiple Fourier series summed over rectangles the Weyl multiplier is the sequence
$$
W(\nu)=\prod_{j=1}^{N}\log(|\nu_j|+2),\quad \nu=(\nu_1,\dots,\nu_N)\in {\mathbb Z}^N,\quad N\geq2.
$$
It is difficult to determine the authorship of this result. Since in the multidimensional case it follows from more general (thoroughly proved) results by F.\,Moricz \cite{Moricz} (1981), so usually F.\,Moricz is considered to be its author. For $N=2$ it was obtained by S.\,Kaczmarz \cite{Kaczmarz} (1930); however, proof of several estimates in \cite{Kaczmarz} causes questions -- the matter concerns the proofs of Lemma 1 on asymptotic of partial sums (p. 93) and of the Theorem (p. 95). Remarks  concerning this see in \cite{Moricz}, \cite{Bloshanskii}. Proof of this result for $N>2$ was made in 1977 by M.\,Kojima \cite[Theorem 3]{Kojima}, but the lemma on asymptotic of partial sums in his paper is given without the proof with the reference that it can be proved the same as in \cite{Kaczmarz} for $N=2$. Note that this result for multidimensional case was stated in 1973 by L.\,V.\,Zhizhiashvili \cite[p. 90]{Gigia}; in 1977 J.\,Chen, N.\,Shieh \cite{Chen} actually stated this result once more (without reference to \cite{Gigia}), remaining the basic estimates in their paper without proofs.
For $N=1$ all the results listed above are analogs of the classical theorem by A.\,N.~Kolmogorov, G.\,A.~Seliverstov, \cite{Kolmogorov2} and A.\,I.\,Plessner \cite{Plessner} (1925-1926).

In the case $N=2$  P.\,Sjolin \cite[Theorem 7.2]{Sjolin} (1971) proved that the following sequence can be taken as the Weyl multiplier
\begin{equation}\label{0.2}
W(\nu_1,\nu_2)=\log^2\bigl[\,\min(|\nu_1|,|\nu_2|)+2\bigl],
\end {equation}
and E.\,M.\,Nikishin \cite[Theorem 4]{Nikishin} (1972) proved that \eqref{0.2} is the exact Weyl multiplier.

{\bf 3.} Let $N\ge1,$ $M=\{1,\dots,N\}$ and $s\in M$. Denote:
$J_s=\{j_1,...,j_s\}$, $j_q<j_l$ for $q<l$, and (in the case $s<N$)
$M\setminus J_s=\{\alpha_1,\dots,\alpha_{N-s}\}$, $\alpha_q<\alpha_l$ for $q<l$,
these are nonempty subsets of the set $M$. We also consider $J_0= M\setminus J_N=\emptyset$.

Fix an arbitrary $k$, $1\leq k< N$, $N\geq2$, and consider a sample $J_k \subset M$. Define the vectors
$$
\lambda = \lambda[J_k] = (\lambda_{j_1},\dots,\lambda_{j_k}) \in  {\mathbb Z}^k_1,\quad j_s \in J_k,\quad s=1,\dots,k,
$$
and
$$
m=m[J_k]=(n_{\alpha_1}, \dots, n_{\alpha_{N-k}}) \in  {\mathbb Z}^{N-k},\quad \alpha_s \in M\setminus J_k,\quad s=1,\dots,N-k.
$$
We will denote by the symbol
$$
n^{(\lambda, m)}=n^{(\lambda, m)}[J_k]=(n_1,\dots,n_N)\in {\mathbb Z}^N_1
$$
such $N$-dimensional vector, whose components $n_j$ with indices $j\in J_k$ are elements of some (single) lacunary sequences,
i.e., for $j \in J_k\colon n_j = n_j^{(\lambda_j)}\in {\mathbb Z}_1^1$, $n_j^{(1)}=1$, $\frac{n_j^{(\lambda_j+1)}}{n_j^{(\lambda_j)}} \ge q_j>1$, $\lambda_j=1, 2,\dots$, and $n_j^{(\lambda_j)}\to
\infty $ as $\lambda_j \to \infty $; we set
\begin{equation}\label{0.3}
q = q(J_k) = (q_{j_1},\dots,q_{j_k})
\in  {\mathbb R}^k,\quad j_s \in J_k,\quad s=1,\dots,k.
\end{equation}
In its turn, the components of the vector $m[J_k]$ are free. Further in the paper a sequence of partial sums of the type $S_{n^{(\lambda, m)}[J_k]}(x;f)$ we will call a "$J_k$- lacunary"\, sequence of partial sums of multiple Fourier series.

Denote
\begin{equation}\label{0.4}
W(m[J_k])=\prod \limits_{j=1}^{N-k} \log(|n_{{\alpha}_j}|+2).
\end {equation}

{\bf Theorem 1.} {\it Let $J_k$ be an arbitrary sample from $M$,
$1 \leq k \leq N-2$, $N\geq3$. For any function $f\in L_2(\mathbb T^N)$
\begin{equation}\label{0.5}
\biggl\|\,\,\sup\limits_{\lambda_j >\,0\, , j \in J_k, \atop {n_j>0, j \in M\setminus J_k}}\frac{|S_{n^{(\lambda, m)}[J_k]}(x;f)|}{\sqrt{W(m[J_k])}}
\biggl\|_{L_2({\mathbb T}^N)}
\leq C\|f\|_{L_2({\mathbb
T}^N)},
\end {equation}
where the constant $C$ does not depend on the function $f$, $C=C(J_k,q)$\,\,\footnote{\,Further we will denote as $C$ the constants, which are, generally speaking, different.} , and the quantity $q$ is defined in \eqref{0.3}.}

{\bf Remark 1.} In the case $k= N-1$ (i.e., one component is free and the rest $N-1$ are lacunary), M. Kojima proved \cite[Theorem 2]{Kojima} that for any function $f \in L_p({\mathbb T}^N), p>1, N\geq 2$, the following estimate is true:
\begin{equation}\label{0.6}
\biggl\|\sup\limits_{\lambda_j >\,0\, , j \in J_{N-1}, \atop {n_j>0, j \in M\setminus J_{N-1}}}{|S_{n^{(\lambda, m)}[J_{N-1}]}(x;f)|}
\biggl\|_{L_p({\mathbb T}^N)} \leq C\|f\|_{L_p({\mathbb T}^N)}.
\end {equation}

{\bf Theorem 2.} {\it Let $J_{k}$ be an arbitrary sample from $M$, $1\leq k\le N-2$, $N\geq3$. If the Fourier coefficients $c_n(f), n\in {\mathbb Z}^N$, of the function $f\in L_2(\mathbb T^N)$ satisfy condition
\begin{equation}\label{0.7}
\Sigma[f,J_k]=\sum_{n\in {\mathbb Z}^N}|c_n(f)|^2\, W(m[J_k])<+\infty,
\end {equation}
then
\begin{equation}\label{0.8}
\lim_{\lambda_j \to \infty,\, j \in J_{k}, \atop {n_j \to \infty,\, j\in M\setminus J_{k}}}S_{n^{(\lambda, m)}[J_{k}]}(x;f)= f(x)
\quad \text {almost everywhere on}\quad \mathbb T^N.
\end {equation}
Moreover, for any $\alpha >0$ the inequality is true
\begin{equation}\label{0.9}
\mu \biggl\{ x \in {\mathbb T}^N: \sup\limits_{\lambda_j >\,0\, , j \in J_k, \atop {n_j>0, j \in M\setminus J_k}}|S_{n^{(\lambda, m)}[J_k]}(x;f)|>\alpha\biggl\} \leq \frac{C}{\alpha ^2}\cdot\Sigma[f,J_k],
\end {equation}
where $\mu$ is the $N$-dimensional Lebesgue measure, and the constant $C$ does not depend on the function  $f$.}

Theorem 2 can be strengthened for $k=N-2$.

{\bf Theorem 3.} {\it Let $J_{N-2}$ be an arbitrary sample from $M$, $N\geq3$. If the Fourier coefficients $c_n(f)$ of the function $f\in L_2(\mathbb T^N)$ satisfy condition
\begin{equation}\label{0.10}
\Sigma_0[f,J_{N-2}]=\sum_{n\in {\mathbb Z}^N}|c_n(f)|^2\log^2\Bigl[\,\min\limits_{i,\, j\in M\setminus J_{N-2}}(|n_i|,|n_j|)+2\Bigr]<+\infty,
\end {equation}
then
\begin{equation}\label{0.11}
\lim_{\lambda_j \to \infty,\, j \in J_{N-2}, \atop {n_j \to \infty,\, j\in M\setminus J_{N-2}}}S_{n^{(\lambda, m)}[J_{N-2}]}(x;f)= f(x)
\quad \text {almost everywhere on}\quad \mathbb T^N;
\end {equation}
moreover,
\begin{equation}\label{0.12}
\biggl\|\,\sup\limits_{\lambda_j >\,0\, , j \in J_{N-2}, \atop {n_j>0, j \in M\setminus J_{N-2}}}|S_{n^{(\lambda, m)}[J_{N-2}]}(x;f)|
\biggl\|^2_{L_2({\mathbb T}^N)}
\leq C\cdot\Sigma_0[f,J_{N-2}],
\end {equation}
where the constant $C$ does not depend on the function $f$.}

\section{Proof of Theorem 1}\label{s2}

Proof of the theorem is based on the ideas represented by us in \cite{Bloshanskii1}; furthermore, for simplicity of understanding of the proof of this theorem, we'll use the structure and notations elaborated by us in the proofs of Lemma 1 and Theorem 1 in \cite{Bloshanskii1}.

In order to prove the theorem it is necessary to prove the following lemma.

{\bf Lemma 1.} {\it Let $J_1=\{r\}, 1\leq r \leq N$. Then for any function $f\in L_2(\mathbb T^{N})$, $N\geq3$,
\begin{equation}\label{1.1}
\biggl\|\sup\limits_{\lambda_r >\,0\, , r \in J_1, \atop {n_j>0, j \in M\setminus J_1}} \frac{|S_{n^{(\lambda, m)}[J_1]}(x;f)|}{\sqrt{W(m[J_1])}}
\biggl\|_{L_2({\mathbb T}^{N})}
\leq C\|f\|_{L_2({\mathbb
T}^{N})},
\end {equation}
where the constant $C$ does not depend on the function  $f$, $C=C(J_1,q)$, and the quantity $q$ is defined in \eqref{0.3}.}

{\bf Proof of Lemma 1.}\footnote{In the proof of this lemma, some ideas represented by P.Sjolin \cite{Sjolin} and M.Kojima \cite{Kojima} are used.}
 Not to complicate the proof, let us consider $r=1$. We denote $\widetilde{x}=(x_2,x_3,\dots,x_N)\in\mathbb T^{N-1},$ and consider
\begin{equation}\label{1.2}
\widetilde{\mathbb T}^{N-1}=\bigl\{\widetilde{x}\in\mathbb T^{N-1}: g(x_1)=f(x_1,\widetilde{x})\in L_2(\mathbb T^1) \bigl\};
\end {equation}
it is obvious,
\begin{equation}\label{1.3}
\mu_{N-1}\widetilde{\mathbb T}^{N-1}=\mu_{N-1}\mathbb T^{N-1}=(2\pi)^{N-1},
\end {equation}
here $\mu_{N-1}$ is the $(N-1)$-dimensional Lebesgue measure.

Fix an arbitrary point $\widetilde{x}\in \widetilde{\mathbb T}^{N-1}$ and expand the function $g(x_1)$ in the (single)
trigonometric Fourier series
\begin{equation}\label{1.4}
g(x_1)\thicksim\sum_{k\in {\mathbb Z}^1}c_ke^{ikx_1}.
\end {equation}
Consider the partial sums of this series $S_m(x_1;g)$ with the indices $m=n_1^{(\lambda_1)}\in {\mathbb Z}_1^1$, $\lambda_1=1,2,\dots$, where $\bigl\{n_1^{(\lambda_1)}\bigl\}$ is a lacunary sequence; set $n_1^{(0)}=0$
and define the difference
$$
\Delta_{\lambda_1}(x_1;g)=\biggl\{\begin{array}{lcl}
  S_0(x_1;g) & \text{for} & \lambda_1=0,\\
  S_{n_1^{(\lambda_1)}}(x_1;g)-S_{n_1^{(\lambda_1-1)}}(x_1;g) & \text{for} &  \lambda_1=1,2,\dots. \\
\end{array} \biggl.
$$

Let us split the series in \eqref{1.4} into two series:
\begin{equation}\label{1.5}
\sum\limits_{\lambda_1=0}^{\infty}\Delta_{2\lambda_1+1}(x_1;g), \quad \sum\limits_{\lambda_1=0}^{\infty}\Delta_{2\lambda_1}(x_1;g).
\end {equation}

From \cite[Ch. 15, Theorem (4.11)]{Zigmund2} it follows that trigonometric series \eqref{1.5} are Fourier series of some functions $g_1(x_1)=f_1(x_1,\widetilde{x})$ and  $g_2(x_1)=f_2(x_1,\widetilde{x})$, $g_1, g_2\in L_2(\mathbb T^1)$ (here we took account of notation \eqref{1.2}), and inequalities are true
\begin{equation}\label{1.6}
\|g_1\|_{L_2({\mathbb T}^1)}
\leq C\|g\|_{L_2({\mathbb
T}^1)}, \quad \|g_2\|_{L_2({\mathbb T}^1)}
\leq C\|g\|_{L_2({\mathbb T}^1)}.
\end {equation}
In its turn, taking into account L.\,Carleson's result \cite{Carleson} (for the one-dimensional trigonometric Fourier series), we have:
$$
g_1(x_1)=\sum\limits_{\lambda_1=0}^{\infty}\Delta_{2\lambda_1+1}(x_1;g), \quad g_2(x_1)=\sum\limits_{\lambda_1=0}^{\infty}\Delta_{2\lambda_1}(x_1;g) \quad \text{for a.e.} \quad x_1\in\mathbb T^1.
$$
Hence, in view of the definition of the functions $g$, $g_1$ and $g_2$ (as well as notation \eqref{1.2}) we obtain
\begin{equation}\label{1.7}
f(x_1,\widetilde{x})=g(x_1)=g_1(x_1)+g_2(x_1)=f_1(x_1,\widetilde{x})+f_2(x_1,\widetilde{x})\quad \text{for a.e.}\quad x_1\in\mathbb T^1.
\end {equation}
In its turn, taking into account that, according to the assumption of the lemma, $f\in L_2(\mathbb T^N)$, in view of estimates \eqref{1.3}, \eqref{1.6} and arbitrariness of the choice of $\widetilde{x}\in \widetilde{\mathbb T}^{N-1}$, we obtain the following estimates
\begin{equation}\label{1.8}
\|f_j \|_{L_2({\mathbb T}^{N})}\leq C\|f\|_{L_2({\mathbb T}^{N})}, \quad j=1,2.
\end {equation}

Now, denoting for convenience
$$
b_{m}=\{W(m)\}^{-\frac12},\quad m=m[J_1],
$$
we define the functions $G_{m}(x_1,\widetilde{x})$, $G_{m}^{(1)}(x_1,\widetilde{x})$ and $G_{m}^{(2)}(x_1,\widetilde{x})$, as follows
\begin{equation}\label{1.9}
G_{m}(x_1,\widetilde{x})=S_{m}(\widetilde{x};f(x_1,\cdot))\cdot b_{m},\quad G_{m}^{(j)}(x_1,\widetilde{x})=S_{m}(\widetilde{x};f_j(x_1,\cdot))\cdot b_{m},\,j=1,2.
\end {equation}
From equality \eqref{1.7} we get:
\begin{equation}\label{1.10}
S_{(n_1^{(\lambda_1)},\, m)}(x;f)\,b_{m}=
S_{n_1^{(\lambda_1)}}(x_1;G_{m}^{(1)}(\cdot,\widetilde{x}))+S_{n_1^{(\lambda_1)}}(x_1;
G_{m}^{(2)}(\cdot,\widetilde{x})).
\end {equation}

Note, that in view of the definition of the functions $f_j(x_1,\widetilde{x})$, $j=1,2,$ in \eqref{1.7}, for any fixed $\widetilde{x}\in \widetilde{\mathbb T}^{N-1}$ the Fourier coefficients of the function $f_1(x_1,\widetilde{x})$ (over the variable $x_1$) $c_k(f_1)=0$ for  $n_1^{(2\lambda_1+1)}<|k|\leq n_1^{(2\lambda_1+2)}$, and the Fourier coefficients of the function $f_2(x_1,\widetilde{x})$ (over the variable $x_1$) $c_k(f_2)=0$ for  $n_1^{(2\lambda_1)}<|k|\leq n_1^{(2\lambda_1+1)}$. In its turn, taking account of the definition of the functions $G_{m}^{(j)}(x_1,\widetilde{x})$, $j=1,2$ (see \eqref{1.9}), the Fourier coefficients of the function $G_{m}^{(1)}(x_1,\widetilde{x})$ (over the variable $x_1$) $c_k(G_{m}^{(1)})=0$ for  $n_1^{(2\lambda_1+1)}<|k|\leq n_1^{(2\lambda_1+2)}$, and the Fourier coefficients of the function $G_{m}^{(2)}(x_1,\widetilde{x})$ (over the variable $x_1$) $c_k(G_{m}^{(2)})=0$ for  $n_1^{(2\lambda_1)}<|k|\leq n_1^{(2\lambda_1+1)}$.

Hence, both functions $G_{m}^{(j)}(x_1,\widetilde{x})$, $j=1,2$ (over the variable $x_1$) satisfy the assumptions of Lemma (1.19) from \cite[Ch. 13]{Zigmund2} (see also \cite[Ch. VI, p. 73]{Hardy}, \cite[Ch. III, p. 79]{Zigmund1}; for the detailed formulation of this statement, appropriate for understanding of the proof, see  \cite[Theorem B]{Bloshanskii1}).
So, in view of this lemma, the following estimates hold true:
\begin{equation}\label{1.11}
\sup\limits_{\lambda_1>\,0}|S_{n_1^{(\lambda_1)}}(x_1;G_{m}^{(j)}(\cdot,\widetilde{x}))|\leq C\sup\limits_{n_1>\,0}|\sigma_{n_1}(x_1;G_{m}^{(j)}(\cdot,\widetilde{x}))|,\quad j=1,2,
\end {equation}
where the constant $C=C(q)$ does not depend on $G_{m}^{(j)}, j=1,2$,
and $\sigma_{n}(t;\varphi)$ are the Cezaro means
\begin{equation}\label{1.12}
\sigma_{n}(t;\varphi)=\frac {1}{n+1}\sum\limits_{r=0}^{n}S_r(t;\varphi),\quad t \in {\mathbb T}^1.
\end{equation}
In its turn, for the Cezaro means \eqref{1.12} the estimate is true (see \cite[Ch. 4, Theorem (7.8)]{Zigmund1}):
 \begin{equation}\label{1.13}
\biggl\|\sup\limits_{n>\,0}|\sigma_{n}(t;\varphi)|
\biggl\|_{L_p({\mathbb T}^1)}
\leq C\|\varphi\|_{L_p({\mathbb
T}^1)},\quad \varphi\in L_p(\mathbb T^1), 1<p<\infty,
\end{equation}
where the constant $C$ does not depend on the function $\varphi$.

By \eqref{1.11} and \eqref{1.13} we have:
$$
\biggl\|\sup\limits_{\lambda_1>0,\,n_j>0,\, j\in M\setminus J_1}|S_{(n_1^{(\lambda_1)}, m)} (x_1;G_{m}^{(j)}(\cdot,\widetilde{x}))|
\biggl\|_{L_2({\mathbb T}^{N})}
$$
$$
\leq C\biggl\|\sup\limits_{n_1>\,0}|\sigma_{n_1}(x_1;\sup\limits_{n_j>0,\, j\in M\setminus J_1} |S_{m}(\widetilde{x};f_j(x_1,\cdot))|\,b_{m})|
\biggl\|_{L_2({\mathbb T}^{N})}
$$
\begin{equation}\label{1.14}
\leq C\biggl\|\sup\limits_{n_j>0,\, j\in M\setminus J_1}|S_{m}(\widetilde{x};f_j(x_1,\cdot))|\cdot b_{m}
\biggl\|_{L_2({\mathbb T}^{N})},\quad j=1,2.
\end {equation}
In view of the result by F. Moricz  \cite [Theorem 1 ]{Moricz} for $\varphi\in L_2({\mathbb T}^{d})$,
$d \geq2$, $n\in {\mathbb Z}_1^{d}$, the estimate is true:
\begin{equation}\label{1.15}
\biggl\|\,\sup\limits_{n_j>\,0,\atop{ j=1,\dots,d}}\frac{|S_n(t;\varphi)|}{\sqrt{W(n)}}\,
\biggl\|_{L_2({\mathbb T}^{d})}
\leq C\|\varphi\|_{L_2({\mathbb
T}^{d})},  \quad W(n)=\prod \limits_{j=1}^d \log(n_j+2).
\end {equation}
From \eqref{1.14} and \eqref{1.15}, considering \eqref{1.9}, we get:
\begin{equation}\label{1.16}
\biggl\|\sup\limits_{\lambda_1>0,\,n_j>0,\, j\in M\setminus J_1}|S_{n_1^{(\lambda_1)}}(x_1;G_{m}^{(j)}(\cdot,\widetilde{x}))|
\biggl\|_{L_2({\mathbb T}^{N})}\leq C(J_1)\|f_j\|_{L_2({\mathbb
T}^{N})},\quad j=1,2.
\end {equation}

Further, from equality \eqref{1.10} and estimates \eqref{1.8} and \eqref{1.16} it follows:
$$
\biggl\|\sup\limits_{\lambda_1>0,\,n_j>0,\, j\in M\setminus J_1}|S_{(n_1^{(\lambda_1)}, m)}(x;f)|\cdot b_{m}
\biggl\|_{L_2({\mathbb T}^{N})}
\leq C(J_1,q)\|f\|_{L_2({\mathbb
T}^{N})}.
$$
Thus, taking account of our assumptions, we prove estimate \eqref{1.1}.

Lemma 1 is proved.

{\bf Proof of Theorem 1.} The proof of estimate \eqref{0.5} will be conducted by the induction on $N$, $N\geq3$.

The first step of induction, i.e. $N=3$; in this case we must prove that
for any $J_1=\{r\}$, $1\leq r\leq 3$, for any function $f\in L_2(\mathbb T^3)$:
\begin{equation}\label{1.17}
\biggl\|\sup\limits_{\lambda_r >\,0\, , \atop {n_j>0,\, j \in M\setminus J_1}}\frac{|S_{n^{(\lambda, m)}[J_1]}(x;f)|}{\sqrt{W(m[J_1])}}
\biggl\|_{L_2({\mathbb T}^3)}
\leq C(J_1,q)\|f\|_{L_2({\mathbb
T}^3)}.
\end {equation}
As we see, the validity of \eqref{1.17} follows from the validity of Lemma 1, i.e., from estimate \eqref{1.1} for $N=3$.

Further, suppose that \eqref{0.5} is true for some $N=l$, $l\geq3$, i.e., for any
$J_k$ in $M=\{1,\dots,l\}$, $1\leq k\leq l-2$, for any function $f\in L_2(\mathbb T^{l})$,
\begin{equation}\label{1.18}
\biggl\|\sup\limits_{\lambda_j >\,0,\, j \in J_k, \atop {n_j>0, \,j \in M\setminus J_k\,} }\frac {|S_{n^{(\lambda,m)}[J_k]}(x;f)|}{\sqrt{W(m[J_k])}}
\biggl\|_{L_2({\mathbb T}^{l})}
\leq C(J_k,q)\|f\|_{L_2({\mathbb T}^{l})}.
\end {equation}

Let us prove that estimate \eqref{0.5} is true for $N=l+1$, i.e., for any $J_d$ in $M=\{1,\dots,l+1\}$, $1\leq d\leq l-1$,
for any function $f\in L_2(\mathbb T^{l+1})$,
\begin{equation}\label{1.19}
\biggl\|\sup\limits_{\lambda_j >\,0,\, j \in J_d, \atop {n_j>0, \,j \in M\setminus J_d}}\frac {|S_{n^{(\lambda, m)}[J_d]}(x;f)|}{\sqrt{W(m[J_d])}}
\biggl\|_{L_2({\mathbb T}^{l+1})}\leq C(J_d, q)\|f\|_{L_2({\mathbb T}^{l+1})}.
\end {equation}

If $d=1$, then \eqref{1.19} again follows from the result of Lemma 1.

Consider now $d\geq 2$, and, to simplify the notation, let us assume that the sample $J_d$ is of the form: $J_d=\{1,2,\dots, d\}
$. In this case,
$$
n^{(\lambda, m)}[J_d]=(n_1^{(\lambda_1)},n_2^{(\lambda_2)},\dots, n_d^{(\lambda_d)},
n_{d+1},\dots,n_{l+1}),\quad
m[J_d]=(n_{d+1},\dots,n_{l+1}).
$$
We denote
$$
\widetilde{n}^{(\lambda, m)}= \widetilde{n}^{(\lambda, m)}[J_d]=(n_2^{(\lambda_2)},\dots,n_d^{(\lambda_d)}, n_{d+1},\dots,n_{l+1})\in{\mathbb Z}_0^{l}.
$$

Let the set $\widetilde{\mathbb T}^{l}$ be defined analogously to \eqref{1.2}, precisely,
\begin{equation}\label{1.20}
\widetilde{\mathbb T}^l=\bigl\{\widetilde{x}=(x_2,x_3,\dots,x_{l+1})\in\mathbb T^l: g(x_1)=f(x_1,\widetilde{x})\in L_2(\mathbb T^l) \bigl\};
\end {equation}
it is obvious, $\mu_l\widetilde{\mathbb T}^l=\mu_l\mathbb T^l=(2\pi)^l$ (here $\mu_l$ is the $l$-dimensional Lebesgue measure).
Fixing an arbitrary point $\widetilde{x}\in \widetilde{\mathbb T}^l$, by the same argumentation as in Lemma 1 (see \eqref{1.5} -- \eqref{1.7}), we define two
functions  $g_1(x_1)=f_1(x_1,\widetilde{x})$ and  $g_2(x_1)=f_2(x_1,\widetilde{x})$, $g_1, g_2\in L_2(\mathbb T^1)$,
\begin{equation}\label{1.21}
f(x_1,\widetilde{x})=g(x_1)=g_1(x_1)+g_2(x_1)=f_1(x_1,\widetilde{x})+f_2(x_1,\widetilde{x})\quad \text{for a.e.}\quad x_1\in\mathbb T^1,
\end {equation}
which satisfy (in account of \cite[Ch. 15, Theorem (4.11)]{Zigmund2}) the estimates:
\begin{equation}\label{1.22}
\|f_j \|_{L_2({\mathbb T}^{l+1})}\leq C\|f\|_{L_2({\mathbb T}^{l+1})}, \quad \quad j=1,2.
\end {equation}

Further, analogously to \eqref{1.9} we define the following functions:
$$
G_{\widetilde{n}^{(\lambda, m)}}(x_1,\widetilde{x})=S_{\widetilde{n}^{(\lambda, m)}}
(\widetilde{x};f(x_1,\cdot))\cdot b_m, \quad G_{\widetilde{n}^{(\lambda, m)}}^{(j)}(x_1,\widetilde{x})=
$$
$$
=S_{\widetilde n ^{(\lambda, m)}}(\widetilde{x};f_j(x_1,\cdot))\cdot b_m,\quad j=1,2,\quad  m=m[J_d],\quad b_m=\{W(m[J_d])\}^{-\frac12}.
$$
From equality \eqref{1.21} we have:
$$
S_{n^{(\lambda, m)}[J_d]}(x;f)\cdot b_m=S_{n_1^{(\lambda_1)}}(x_1;G_{\widetilde{n}^{(\lambda, m)}}(\cdot,\widetilde{x}))
$$
\begin{equation}\label{1.23}
=S_{n_1^{(\lambda_1)}}(x_1;G_{\widetilde{n}^{(\lambda, m)}}^{(1)}(\cdot,\widetilde{x}))+
S_{n_1^{(\lambda_1)}}(x_1;G_{\widetilde{n}^{(\lambda, m)}}^{(2)}(\cdot,\widetilde{x})).
\end {equation}
By the same argumentation as in the proof of \eqref{1.11}, we obtain:
$$
\sup\limits_{\lambda_1>\,0}|S_{n_1^{(\lambda_1)}}(x_1;G_{\widetilde{n}^{(\lambda, m)}}^{(j)}(\cdot,\widetilde{x}))|\leq C\sup\limits_{n_1>\,0}|\sigma_{n_1}(x_1;G_{\widetilde{n}^{(\lambda, m)}}
^{(j)}(\cdot,\widetilde{x}))|,\quad j=1,2.
$$
The same as in the proof of Lemma 1, using estimate \eqref{1.13}, we get:
$$
\biggl\|\sup\limits_{\lambda_1,\,\dots,\lambda_{d}\,>0, \atop{n_j>\,0, j \in M\setminus J_d}}
|S_{n_1^{(\lambda_1)}}(x_1;G_{\widetilde{n}^{(\lambda, m)}}^{(j)}(\cdot,\widetilde{x}))|
\biggl\|_{L_2({\mathbb T}^{l+1})}
$$
$$
\leq C\biggl\|\sup\limits_{n_1>\,0}|\sigma_{n_1}(x_1;\sup\limits_{\lambda_2,\,\dots,\lambda_{d}\,>0,
 \atop{n_j>\,0, j \in M\setminus J_d}}
|S_{\widetilde{n}^{(\lambda, m)}}(\widetilde{x};f_j(x_1,\cdot))|\cdot b_{m})
\biggl\|_{L_2({\mathbb T}^{l+1})}
$$
\begin{equation}\label{1.24}
\leq C\biggl\|\sup\limits_{\lambda_2,\,\dots,\lambda_{d}>0,\atop{n_j>\,0, j \in M\setminus J_d}}
|S_{n_2^{(\lambda_2)},\dots,n_{d}^{(\lambda_{d})},\,n_{d+1},\dots,\,n_{l+1}}(\widetilde{x};f_j(x_1,\cdot))|
\cdot b_{m}
\biggl\|_{L_2({\mathbb T}^{l+1})},\ j=1,2.
\end {equation}

Note that $\{n_j^{(\lambda_j)}\}$, $n_j^{(\lambda_j)}\in {\mathbb Z}_0^1$,
$\lambda_j=1,2,\dots$, $j=2,\dots,d$, are lacunary sequences, and also $1\leq d-1\leq l-2$, and the functions $f_j(x_1,\widetilde{x})\in L_2({\mathbb T}^{l+1})$, $j=1,2$. So, in order to estimate the right part of \eqref{1.24}, we can use the inductive proposition, i.e., the majorant estimate \eqref{1.18}, namely:
$$
\biggl\|\sup\limits_{\lambda_2,\,\dots,\lambda_{d}>\,0, \atop{n_j>\,0, j \in M\setminus J_d}}
|S_{n_2^{(\lambda_2)},\dots,n_{d}^{(\lambda_{d})},\,n_{d+1},\dots,\,n_{l+1}}
(\widetilde{x};f_j(x_1,\cdot))|\cdot b_m \biggl\|^2_{L_2({\mathbb T}^{l+1})}
$$
$$
=\int\limits_{{\mathbb T}^1}\biggl\{ \int\limits_{{\mathbb T}^l}\biggl\{\sup\limits_{\lambda_2,\,\dots,\lambda_{d}>\,0, \atop{n_j>\,0, j \in M\setminus J_d}}
|S_{n_2^{(\lambda_2)},\dots,n_{d}^{(\lambda_{d})},\,n_{d+1},\dots,\,n_{l+1}}
(\widetilde{x};f_j(x_1,\cdot))|\cdot b_m \biggl\}^2 d{\widetilde{x}}\biggl\}dx_1
$$
$$
\leq C \int\limits_{{\mathbb T}^1}\biggl\{\int\limits_{{\mathbb T}^l}|f_j(x_1,\widetilde{x})|^2 d{\widetilde{x}}\biggl\}dx_1 = C\|f_j\|^2_{L_2({\mathbb T}^{l+1})},\quad j=1,2.
$$
By this and \eqref{1.24} we have:
\begin{equation}\label{1.25}
\biggl\|\sup\limits_{\lambda_1,\,\dots,\lambda_{d}>\,0, \atop{n_j>\,0, j \in M\setminus J_d}}
|S_{n_1^{(\lambda_1)}}(x_1;G_{\widetilde{n}^{(\lambda,m)}}^{(j)}(\cdot,\widetilde{x}))|
\biggl\|_{L_2({\mathbb T}^{l+1})}\leq C\|f_j\|_{L_2({\mathbb
T}^{l+1})},\quad j=1,2.
\end {equation}
Further, from \eqref{1.23}, \eqref{1.22} and \eqref{1.25} it follows the validity of estimate \eqref{1.19}.

In view of the induction method, we get that estimate \eqref{0.5} is true for any $N\geq3$ and any $k$ (the number of lacunary components in the vector $n^{(\lambda,m)}[J_k]\in{\mathbb Z}_0^N$), $1\leq k\leq N-2$.

Theorem 1 is proved.

A simple corollary of Theorem 1 is the following statement which will be used in the proof of Theorem 2.

Let $J_k$ be an arbitrary sample from $M$, $1 \leq k \leq N-2$, $N\geq3$. Fix an integer $s, 1\le s \le N-k,$ and indices $\nu_1,\dots,\nu_s \in M\setminus J_k$.
Denote
\begin{equation}\label{1.26}
Q_{n^{(\lambda, m)}[J_k]}^{(p_{\nu_1},\dots,p_{\nu_s})}(x;f)= \sum \limits_{n_{\nu_1}=0}^{p_{\nu_1}} \dots \sum \limits_{n_{\nu_s}=0}^{p_{\nu_s}} S_{n^{(\lambda, m)}[J_k]}(x;f),\quad p_{\nu_1},\dots,p_{\nu_s} \in {\mathbb Z}_1^1.
\end {equation}

{\bf Corollary of Theorem 1.} {\it For any function $f\in L_2(\mathbb T^N)$, for any $s, 1\le s \le N-k$, for any $\nu_1,\dots,\nu_s \in M\setminus J_k$ and $p_{\nu_1},\dots,p_{\nu_s}\in {\mathbb Z}_1^1,$ the following estimate is true
\begin{equation}\label{1.27}
Q_{n^{(\lambda, m)}[J_k]}^{(p_{\nu_1},\dots,p_{\nu_s} )}(x;f)= O\biggl(\prod\limits_{j=1}^s \, p_{\nu_j} \prod\limits_{l \in M\setminus J_k ,\atop{l\neq \nu_1,\dots,\nu_s}}\sqrt{\log(|n_{l}|+2)}\biggl)\quad \text{for a.e.}\quad x\in \mathbb T^N.
\end {equation}}

{\bf Proof of Corollary of Theorem 1}. Let us prove estimate \eqref{1.27} for $s=1$ (for $s>1$ the proof is similar). Let for definiteness $J_k=\{1,\dots,k\}$ and $\nu_1=k+1$, then $n^{(\lambda, m)}[J_k]=(n_1^{(\lambda_1)},\dots,n_k^{(\lambda_k)}, n_{k+1},\dots,n_N).$ Set $p=p_{k+1}$ and consider
$$
Q_{n^{(\lambda, m)}[J_k]}^{(p)}(x;f)= \sum \limits_{n_{k+1}=0}^{p}S_{n_1^{(\lambda_1)},\dots,n_k^{(\lambda_k)}, n_{k+1},\dots,n_N}(x;f)
$$
$$
=(p+1)\pi^{-1}\int\limits_{\mathbb T^1}K_{p}(u_{k+1})\biggl\{{\pi^{-N+1}}\int\limits_{\mathbb T^{N-1}}f(u+x)\prod \limits_{j=1}^k D_{n_j^{(\lambda_j)}}(u_j)
$$
\begin{equation}\label{1.28}
\times\prod \limits_{j=k+2}^N D_{n_j}(u_j)d\tilde{u}\biggl\}du_k, \quad \tilde{u}=(u_1,\dots,u_k,u_{k+2},\dots,u_N),
\end{equation}
where $D_{n_j}(u_j)$ is the Dirichlet kernel, $K_p(u_{k+1})$ is the Fej\'{e}r kernel.

Denote as  $G(x_1,\dots,x_k,x_{k+2},\dots,x_N; u_{k+1})$ the expression in the braces in \eqref{1.28}.
Considering notation \eqref{1.12}
$Q_{n^{(\lambda, m)}[J_k]}^{(p)}(x;f)= (p+1)\,\sigma_p (x;G)$.
Hence, for $\sigma_p (x;G)$ estimate \eqref{1.13} is true. Further, note that $G$ is the partial sum of the Fourier series of the function $f$ whose index has $k$ lacunary components, and thus, for $G$ the estimate from Theorem 1 is true. Thus, \eqref{1.27} is proved.

Corollary of Theorem 1 is proved.

\section{Proof of Theorem 2}\label{s2}

In the proof of the theorem ideas from the papers \cite{Kaczmarz} and \cite[Theorem 3]{Kojima} are used.

Let us fix a sample $J_k$, $1\le k \le N-2$. Without loss of generality, let us consider that $J_k=\{N-k+1,\dots,N\}$. In this case the vector $m[J_k]=(n_{1},\dots, n_{N-k})$. Consider
$$
n^{(\lambda,m)}[J_k]=(n_1,\dots,n_{N-k}, n_{N-k+1}^{(\lambda_{N-k+1})},\dots,n_N^{(\lambda_N)}) = ( n_{1},\dots,n_{N-k}, n'),
$$
where
\begin{equation}\label{2.1}
n'=(n_{N-k+1}^{(\lambda_{N-k+1})},\dots,n_N^{(\lambda_N)}).
\end {equation}
Let condition \eqref{0.7} be satisfied, i.e., in our case
\begin{equation}\label{2.2}
\sum_{n\in {\mathbb Z}^N}|c_n(f)|^2\, \prod \limits_{\alpha=1}^{N-k} \log(|n_{\alpha}|+2)<+\infty.
\end {equation}
Estimate \eqref{2.2} permits to assert (see, e.g., \cite[Lemma 3]{Kaczmarz}) that there exists a sequence of numbers $\{p_{j}\},$ $p_{j}>0, p_{j}=p_{-j}$,
increasing to infinity as slowly as we like as $j \to \infty$, such that
\begin{equation}\label{2.3}
\sum_{n\in {\mathbb Z}^N}|c_n(f)|^2\, \prod \limits_{\alpha =1}^{N-k} \log(|n_{\alpha}|+2)\, p_{n_{\alpha}}<+\infty.
\end {equation}
Set
\begin{equation}\label{2.4}
b_j=\bigl\{ \log(|j|+2)\, p_j\bigl\}^{-\frac12},\quad j\in {\mathbb Z}^1.
\end {equation}
Note that $\{b_j\}$ (taking into account the choice of the sequence $\{p_j\}$) is a convex sequence, satisfying the following conditions (see \cite[Ch. III, p. 93]{Zigmund1})
\begin{equation}\label{2.5}
b_j=b_{-j},\quad b_j \to 0 \quad \text {and}\quad j\Delta^1 b_j \to 0 \quad \text {as}\quad j \to \infty,
\quad \sum \limits_{j=1}^{\infty} j \Delta^2 b_j < \infty,
\end {equation}
where $\Delta^0 b_j=b_j,\quad \Delta^1 b_j=b_j-b_{j+1},\quad \Delta^2 b_j=\Delta^1 b_j - \Delta^1 b_{j+1},\quad j\in {\mathbb Z}^1$.

Let $ M_j \in {\mathbb Z}^1_1$ be such that inequalities are true: $2^{M_j^2}\le n_j \le 2^{(M_j+1)^2}$ (i.e., $M_j^2\sim \log n_j,\quad j=1,\dots,N-k$). Denote $\alpha_j= 2^{M_j^2} $ and represent the $J_k$- lacunary partial sum $S_{n^{(\lambda,m)}[J_{k}]}(x;f)$ as follows:
$$
S_{n^{(\lambda,m)}[J_{k}]}(x;f) = \{S_{n_1,n_2,\dots,n_{N-k}, n'}(x;f) - S_{\alpha_1, n_2,\dots, n_{N-k}, n'}(x;f)\}
$$
$$
+ \sum \limits_{j=2}^{N-k-1}\{S_{\alpha_1,\dots,\alpha_{j-1},n_{j},\dots,n_{N-k}, n'}(x;f)- S_{\alpha_1,\dots,\alpha_{j-1},\alpha_{j},n_{j+1},\dots, n_{N-k}, n'}(x;f)\}
$$
\begin{equation}\label{2.6}
+ S_{\alpha_1,\dots,\alpha_{N-k-1}, n_{N-k}, n'}(x;f)=\sum \limits_{j=1}^{N-k-1} \Delta S_{n_j, \alpha_j}(x;f)+ S_{\alpha_1,\dots,\alpha_{N-k-1}, n_{N-k}, n'}(x;f).
\end {equation}
Because the index $(\alpha_1, \dots, \alpha_{N-k-1}, n_{N-k},n') \in{\mathbb Z}^N$ of the latter partial sum has $N-1$ lacunary components (see \eqref{2.1}), by Kojima's result \cite[Theorem 2]{Kojima} (see Remark 1 in the Introduction) the equality holds true
\begin{equation}\label{2.7}
\lim \limits_{\alpha_1\to \infty, \dots,\alpha_{N-k-1}\to \infty, \atop{n_{N-k} \to \infty, n'\to \infty}} S_{\alpha_1,\dots,\alpha_{N-k-1}, n_{N-k}, n'}(x;f)= f(x) \quad \text{a.e. on} \quad {\mathbb T}^N.
\end {equation}
Thus, theorem will be proved if we prove that each difference in \eqref{2.6} tends to zero a.e. on ${\mathbb T}^N$, precisely,
\begin{equation}\label{2.8}
\lim_{n_l\to \infty, l=1,\dots,N-k, \atop{n' \to \infty}} \Delta S_{n_j, \alpha_j}(x;f)=0 \quad \text{a.e. on} \quad {\mathbb T}^N, \quad j=1,\dots,N-k-1.
\end {equation}
Let us prove \eqref{2.8} for $j=1$, for other differences the proof is similar. So, consider
\begin{equation}\label{2.9}
\Delta S_{n_1, \alpha_1}(x;f) = S_{n_1,n_2,\dots,n_{N-k}, n'}(x;f) - S_{\alpha_1,n_2,\dots,n_{N-k}, n'}(x;f).
\end {equation}

In view of \eqref{2.3} and \eqref{2.4}
\begin{equation}\label{2.10}
\sum_{n\in {\mathbb Z}^N}{|c_n(f)|^2}\prod \limits_{j =1}^{N-k} b_{n_j}^{-2} <+\infty;
\end {equation}
so, by the Riezc-Fischer theorem, there exists a function $g \in L_2({\mathbb T}^N)$ such that
\begin{equation}\label{2.11}
c_n(f) = c_n(g)\prod \limits_{j =1}^{N-k} b_{n_j}.
\end {equation}
Let us represent the $J_k$- lacunary partial sum of the Fourier series of the function $g$ in the real form and denote: $I=(i_1,\dots,i_{N-k})$, $I'=(i_{N-k+1},\dots,i_N)$, $i_j\in {\mathbb Z}^1_0,\quad j=1,\dots,N$; we get
\begin{equation}\label{2.12}
S_{n^{(\lambda,m)}[J_{k}]}(x;g)=\sum\limits_{i_1=0}^{n_1}\dots \sum\limits_{i_{N-k}=0}^{n_{N-k}}\sum\limits_{i_{N-k+1}=0}^{n_{N-k+1}^{(\lambda_{N-k+1})}}\dots \sum\limits_{i_{N}=0}^{n_{N}^{(\lambda_{N})}}G_{I, I'}(x;g)=\sum\limits_{i_1=0}^{n_1}\dots \sum\limits_{i_{N-k}=0}^{n_{N-k}}A_{I},
\end {equation}
where
\begin{equation}\label{2.13}
A_{I}=A_{i_1,\dots, i_{N-k}}=\sum\limits_{i_{N-k+1}=0}^{n_{N-k+1}^{(\lambda_{N-k+1})}}\dots \sum\limits_{i_{N}=0}^{n_{N}^{(\lambda_{N})}} G_{I, I'}(x;g).
\end {equation}

Thus, by \eqref{2.12}, \eqref{2.13}, the $J_k$- lacunary partial sum of the Fourier series of the function $f$ in the real form looks as
\begin{equation}\label{2.14}
S_{n^{(\lambda,m)}[J_{k}]}(x;f) =
\sum\limits_{i_1=0}^{n_1}\dots \sum\limits_{i_{N-k}=0}^{n_{N-k}}A_{I}\cdot b_{i_1}\cdot {\dots} \cdot b_{i_{N-k}} .
\end {equation}
Given \eqref{2.14}, the difference in \eqref{2.9} looks as
$$
\Delta S_{n_1, \alpha_1}(x;f) = \sum\limits_{i_1=0}^{n_1} \sum\limits_{i_2=0}^{n_2}{\dots} \sum\limits_{i_{N-k}=0}^{n_{N-k}}A_{I}\cdot b_{i_1}\cdot {\dots} \cdot b_{i_{N-k}}
$$
\begin{equation}\label{2.15}
-\sum\limits_{i_1=0}^{\alpha_1} \sum\limits_{i_2=0}^{n_2}{\dots} \sum\limits_{i_{N-k}=0}^{n_{N-k}}A_{I}\cdot b_{i_1}\cdot {\dots} \cdot b_{i_{N-k}}=B(n_1)-B(\alpha_1).
\end {equation}

Set $\nu=N-k$ and introduce the following notation.

Let $s,p,q$ be integers, $0\le s,p,q \le \nu$, $s+p+q \le \nu$. Denote
$$
{\mathfrak A}(s,p,q)=\{L=(l_1,\dots,l_{\nu})\in {\mathbb Z}^{\nu}_0: 0=l_0<l_1<\dots<l_{s}\le \nu;
$$
$$
0= l_0<l_{s+1}<\dots<l_{s+p}\le \nu;\quad 0=l_0<l_{s+p+1}< \dots< l_{s+p+q} \le \nu;
$$
\begin{equation}\label{2.16}
1\le l_{s+p+q+1}< \dots<l_{\nu}\le\nu; \quad l_{\mu_1}\neq l_{\mu_2} \quad \text {for}\quad {\mu_1}\neq {\mu_2}\}.
\end {equation}
For any vector $ L \in{\mathfrak A}(s,p,q)$ we define the vector
$$
R=R(L)=(r_1,\dots,r_{\nu})\in {\mathbb Z}^{\nu}_0: \quad r_{l_j}=0 \quad \text {for}\quad j=1,\dots, s;
$$
$$
r_{l_j}=1 \quad \text {for}\quad j = s+1,\dots,s+p;\quad r_{l_j}=2\quad\text {for}\quad j=s+p+1,\dots,s+p+q;
$$
\begin{equation}\label{2.17}
r_{l_j}=3 \quad \text {for}\quad j=s+p+q+1,\dots,\nu.
\end {equation}
For any vector $ \varkappa=(\varkappa_{s+1},\dots, \varkappa_{\nu})\in {\mathbb Z}^{\nu-s}_0$ we define the vector
$$
I_\varkappa= (i_1,\dots, i_{\nu}): \quad i_{l_j}=\varkappa_{j} \quad \text {for}\quad j = s+1,\dots,\nu, \quad \text {in the case}\quad 0\le s<\nu;
$$
and set $I_0\equiv I=(i_1,\dots, i_{\nu})$ in the case $s=\nu$.

Let $\{A_{I}\}=\{A_{i_1,\dots,i_{\nu}}\}, \quad A_{i_1,\dots,i_{\nu}}\in {\mathbb R}^1, i_l \in {\mathbb Z}^1_0,\quad l=1,\dots,{\nu},\quad {\nu}\ge 1$, be an arbitrary sequence, and let $\{b_{\alpha}\}, {\alpha} \in {\mathbb Z}^1$, be a convex sequence of real numbers, satisfying conditions \eqref{2.5}. Let elements of the sequence $\{A_{I_{\varkappa}}^{(R)}\}$ be defined as follows:
$$
A_{I_{\varkappa}}^{(R)}= \sum\limits_{i_{l_{s+1}}=0}^{{\varkappa}_{s+1}} \dots \sum\limits_{i_{l_{s+p}}=0}^{{\varkappa}_{s+p}}\,\,
\sum\limits_{t_{s+p+1}=0}^{{\varkappa}_{s+p+1}}\sum\limits_{i_{l_{s+p+1}}=0}^{t_{s+p+1}}{\dots}
\sum\limits_{t_{s+p+q}=0}^{{\varkappa}_{s+p+q}}\sum\limits_{i_{l_{s+p+q}}=0}^{t_{s+p+q}}\,\,
\prod\limits_{j =s+p+q+1}^{\nu} \frac {1}{\Delta^2 b_{{\varkappa}_j}}
$$
\begin{equation}\label{2.18}
\times \sum\limits_{\alpha_{l_{s+p+q+1}}=0}^{{\varkappa}_{s+p+q+1}}\Delta^2 b_{\alpha_{l_{s+p+q+1}}} \sum\limits_{t_{s+p+q+1}=0}^{\alpha_{l_{s+p+q+1}}}
\sum\limits_{i_{l_{s+p+q+1}}=0}^{t_{s+p+q+1}}{\dots}
\sum\limits_{\alpha_{l_{\nu}}=0}^{{\varkappa}_{\nu}}\Delta^2 b_{\alpha_{l_{\nu}}} \sum\limits_{t_{\nu}=0}^{\alpha_{l_{\nu}}}
\sum\limits_{i_{l_{\nu}}=0}^{t_{\nu}}
A_{I}.
\end {equation}
Here we assume that in \eqref{2.18}: there are no sums of the type $ \sum\limits_{i_{l_j}=0}^{\varkappa_j}$ in the case $p=0$; no sums of the type $ \sum\limits_{t_j=0}^{\varkappa_j}\sum\limits_{i_{l_j}=0}^{t_j}$ in the case $q=0$; no sums of the type $ \sum\limits_{\alpha_{l_j}=0}^{\varkappa_j}\Delta^2 b_{\alpha_{l_j}} \sum\limits_{t_j=0}^{\alpha_{l_j}}
\sum\limits_{i_{l_j}=0}^{t_j}$, $j=s+1,\dots,\nu,$ in the case $s+p+q=\nu$.

In particular, by \eqref{2.18} we have: $A_I^{(0,\dots,0)}=A_I$,
$$
A_{\varkappa_1,\dots,\varkappa_{\nu}}^{(1,\dots,1)}=\sum\limits_{i_1=0}^{\varkappa_1}\dots \sum\limits_{i_{\nu}=0}^{\varkappa_{\nu}}A_{i_1,\dots,i_{\nu}},\quad
A_{\varkappa_1,\dots,\varkappa_{\nu}}^{(2,\dots,2)}=\sum\limits_{t_1=0}^{\varkappa_1}\sum\limits_{i_1=0}^{t_1}\dots \sum\limits_{t_{\nu}=0}^{\varkappa_{\nu}}\sum\limits_{i_{\nu}=0}^{t_{\nu}}A_{i_1,\dots,i_{\nu}},
$$
\begin{equation}\label{2.19}
A_{\varkappa_1,\dots,\varkappa_{\nu}}^{(3,\dots,3)}=\prod \limits_{l=1}^{{\nu}} \frac {1}{\Delta^2 b_{\varkappa_l}}\cdot
\sum\limits_{i_1=0}^{\varkappa_1}{\dots}\sum\limits_{i_{\nu}=0}^{\varkappa_{\nu}} A_{i_1,\dots,i_{\nu}}^{(2,\dots,2)}\cdot
\Delta^2 b_{i_1}\cdot {\dots} \cdot\Delta^2 b_{i_{\nu}} .
\end {equation}

{\bf Proposition 1.} {\it For any $n_1,\dots,n_{\nu} \in {\mathbb Z}^1_0$ the following equality holds true
\begin{equation}\label{2.20}
\sum\limits_{i_1=0}^{n_1}\dots \sum\limits_{i_{\nu}=0}^{n_{\nu}} A_{I}\,b_{i_1}\cdot {\dots} \cdot b_{i_{\nu}} =\sum\limits_{r_1=0}^2\dots\sum\limits_{r_{\nu}=0}^2 \Delta^{r_1} b_{n_1-r_1} \cdot {\dots} \cdot \Delta^{r_{\nu}} b_{n_{\nu}-r_{\nu}}A_{n_1-r_1,\dots,n_{\nu}-r_{\nu}}^{(r_1+1,\dots,r_{\nu}+1)}.
\end {equation}
}

Proof of Proposition 1. First, consider $\nu=1$. By setting in \eqref{2.18}, \eqref{2.19} $\nu=1$ and $\varkappa_1=n_1$, we get
\begin{equation}\label{2.21}
A_{n_1}^{(1)}=\sum\limits_{i_1=0}^{n_1} A_{i_1},\quad A_{n_1}^{(2)}=\sum\limits_{t_1=0}^{n_1}\sum\limits_{i_1=0}^{t_1} A_{i_1},\quad
A_{n_1}^{(3)}=\frac {1}{\Delta^2 b_{n_1} }\sum\limits_{i_1=0}^{n_1} A_{i_1}^{(2)}\, \Delta^2 b_{i_1}.
\end {equation}
By applying twice the Abel transformation and considering notation \eqref{2.21} we obtain:
$$
\sum\limits_{i_1=0}^{n_1} A_{i_1} b_{i_1}=\sum\limits_{i_1=0}^{n_1-1} \sum\limits_{k=0}^{i_1}A_k \Delta b_{i_1} + b_{n_1} \sum\limits_{i_1=0}^{n_1}A_{i_1}= \sum\limits_{i_1=0}^{n_1-2} \sum\limits_{l=0}^{i_1} \sum\limits_{k=0}^{l}A_k \Delta^2 b_{i_1}
$$
$$
+\Delta b_{n_1-1}\sum\limits_{i_1=0}^{n_1-1}\sum\limits_{k=0}^{i_1}A_k + b_{n_1}\sum\limits_{i_1=0}^{n_1}A_{i_1} =
\Delta^2 b_{n_1-2} A^{(3)}_{n_1-2} + \Delta b_{n_1-1} A^{(2)}_{n_1-1} + b_{n_1} A^{(1)}_{n_1}.
$$
Finally, we get that estimate \eqref{2.20} is true for $\nu=1$:
\begin{equation}\label{2.22}
\sum\limits_{i_1=0}^{n_1} A_{i_1} b_{i_1}= \sum\limits_{r_1=0}^2 \Delta^{r_1} b_{n_1-r_1} A_{n_1-r_1}^{(r_1+1)}.
\end {equation}
Proposition 1 for $\nu >1$ is proved by application of formula \eqref{2.22} on each index $i_l, \, l=1,\dots,\nu$.

Proposition 1 is proved.

Let us estimate the difference in \eqref{2.15} using formula \eqref{2.20} with $\nu=N-k$; we have
$$
B(n_1)- B(\alpha_1)= \sum\limits_{r_1=0}^2\dots\sum\limits_{r_{\nu}=0}^2 \{Q_{R}(n_1) - Q_{R}(\alpha_1)\}, \quad R=(r_1,\dots,r_{\nu}),
$$
where
\begin{equation}\label{2.23}
Q_{R}(z)=\Delta^{r_1} b_{z-r_1}\cdot \Delta^{r_2} b_{n_2-r_2}\cdot {\dots} \cdot \Delta^{r_{\nu}} b_{n_{\nu}-r_{\nu}}  A_{z-r_1,n_2-r_2,\dots,n_{\nu}-r_{\nu}}^{(r_1+1,r_2+1,\dots,r_{\nu}+1)},\quad z=n_1, \alpha_1.
\end {equation}

We introduce a set
\begin{equation}\label{2.24}
\Omega =\{R(L): L \in \bigcup_{0\,\leq\, s,\,p,\,q\,\leq\,\nu\,, \atop {s+p+q=\nu}}{\mathfrak A}(s,p,q)\}\setminus \{(2,\dots,2)\},
\end {equation}
and by \eqref{2.23},\eqref{2.24}, write the difference \eqref{2.15} as follows
\begin{equation}\label{2.25}
\Delta S_{n_1, \alpha_1}(x;f)= \sum\limits_{R\in \Omega} \{Q_{R}(n_1) - Q_{R}(\alpha_1)\} + \{Q_{(2,\dots,2)}(n_1)-Q_{(2,\dots,2)}(\alpha_1) \}.
\end {equation}

Consider a vector $R\in \Omega$. Without loss of generality, we can consider that $r_1=\dots=r_{s}=0,\,\, r_{s+1}=\dots=r_{s+p}=1, \, r_{s+p+1}=\dots=r_{\nu}=2.$ Denote as $R+1=(r_1+1,\dots,r_{\nu}+1)=(\tilde{r}_1,\dots,\tilde{r}_{\nu})$. Taking into account the choice of the vector $R$, we have:
$\tilde{r}_1=\dots=\tilde{r}_{s}=1,\,\, \tilde{r}_{s+1}=\dots=\tilde{r}_{s+p}=2, \, \tilde{r}_{s+p+1}=\dots=\tilde{r}_{\nu}=3$.
In this case, considering the definition of
$Q_{R}(n_1)$ by \eqref{2.23}, we can write
\begin{equation}\label{2.26}
Q_{R}(n_1) = \prod\limits_{i=1}^{s}b_{n_i}\,\cdot \prod\limits_{i=s+1}^{s+p}\Delta^{1} b_{n_{i}-1}\, \cdot \, \prod\limits_{i=s+p+1}^{\nu}\Delta^{2} b_{n_{i}-2} \, A_{n_1-r_1,\dots,n_{\nu}-r_{\nu}}^{(R+1)}.
\end {equation}

In its turn, taking account of notation \eqref{2.18}, we have:
$$
A_{n_1-r_1,\dots,n_{\nu}-r_{\nu}}^{(R+1)}=
A_{n_1,\dots,n_s,n_{s+1}-1,\dots,n_{s+p}-1,n_{s+p+1}-2,\dots,n_{\nu}-2}^{(R+1)}
$$
$$
=\biggl\{\sum\limits_{i_1=0}^{n_1} \dots \sum\limits_{i_s=0}^{n_s}\quad \biggl \{\sum\limits_{l_{s+1}=0}^{n_{s+1}-1} \sum\limits_{i_{s+1}=0}^{l_{s+1}}\dots \sum\limits_{l_{s+p}=0}^{n_{s+p}-1} \sum\limits_{i_{s+p}=0}^{l_{s+p}}\biggl\{\frac {1}{\Delta^2 b_{n_{s+p+1}-2}}
$$
$$
\times \sum\limits_{\mu_{s+p+1}=0}^{n_{s+p+1}-2} \Delta^2 b_{\mu_{s+p+1}}  \sum\limits_{l_{s+p+1}=0}^{\mu_{s+p+1}} \sum\limits_{i_{s+p+1}=0}^{l_{s+p+1}} \dots
\biggl\{\frac {1}{\Delta^2 b_{n_{\nu}-2}} \sum\limits_{\mu_{\nu}=0}^{n_{\nu}-2} \Delta^2 b_{\mu_{\nu}}
$$
$$
\times\sum\limits_{l_{\nu}=0}^{\mu_{\nu}} \sum\limits_{i_{\nu}=0}^{l_{\nu}} A_{i_1,\dots, i_{\nu}}\biggl\}\biggl\}\biggl\}\biggl\}=\prod\limits_{i=s+p+1}^{\nu}(\Delta^2 b_{n_{i}-2})^{-1}
\sum\limits_{\mu_{s+p+1}=0}^{n_{s+p+1}-2} \Delta^2 b_{\mu_{s+p+1}}\dots
$$
$$
\dots\sum\limits_{\mu_{\nu}=0}^{n_{\nu}-2} \Delta^2 b_{\mu_{\nu}}\cdot\biggl[\sum\limits_{l_{s+1}=0}^{n_{s+1}-1}\dots \sum\limits_{l_{s+p}=0}^{n_{s+p}-1}
\sum\limits_{l_{s+p+1}=0}^{\mu_{s+p+1}}\dots \sum\limits_{l_{\nu}=0}^{\mu_{\nu}}
\biggl\{\sum\limits_{i_1=0}^{n_1} \dots
$$
\begin{equation}\label{2.27}
\dots \sum\limits_{i_s=0}^{n_s} \sum\limits_{i_{s+1}=0}^{l_{s+1}}\dots\sum\limits_{i_{s+p}=0}^{l_{s+p}} \sum\limits_{i_{s+p+1}=0}^{l_{s+p+1}} \dots \sum\limits_{i_{\nu}=0}^{l_{\nu}} A_{i_1,\dots, i_{\nu}}\biggl\}\biggl].
\end {equation}
Note that the expression in the braces equals $S_{n_1,\dots,n_s,l_{s+1},\dots,l_{\nu}, n'}(x;g)$, in view of $\nu=N-k$ (see \eqref{2.12}). Denote the expression in the square brackets by $\Phi=\Phi(\mu_{s+p+1},\dots,\mu_{\nu})$. According to what was said above, we have
$$
\Phi= \sum\limits_{l_{s+1}=0}^{n_{s+1}-1}\dots \sum\limits_{l_{s+p}=0}^{n_{s+p}-1} \sum\limits_{l_{s+p+1}=0}^{\mu_{s+p+1}}\dots \sum\limits_{l_{\nu}=0}^{\mu_{\nu}}S_{n_1,\dots,n_s,l_{s+1},\dots,l_{\nu},n'}(x;g).
$$
By Corollary of Theorem 1 we have
\begin{equation}\label{2.28}
\Phi=O\biggl( \prod\limits_{i=1}^{s}\sqrt{\log n_i}\cdot \prod\limits_{i=s+1}^{s+p} (n_i-1)\biggl)\cdot \mu_{s+p+1}\cdot{\dots}\cdot\mu_{\nu}.
\end {equation}
Thus, \eqref{2.27} and \eqref{2.28} permit us to write
$$
A_{n_1-r_1,\dots,n_{\nu}-r_{\nu}}^{(r_1+1,\dots,r_{\nu}+1)}= \prod\limits_{i=s+p+1}^{\nu}(\Delta^2 b_{n_{i}-2})^{-1}
\cdot O\biggl(\prod\limits_{i=1}^{s}\sqrt{\log n_i}\,\prod\limits_{i=s+1}^{s+p} (n_i-1)\biggl)
$$
\begin{equation}\label{2.29}
\times \sum\limits_{\mu_{s+p+1}=0}^{n_{s+p+1}-2} \mu_{s+p+1}\Delta^2 b_{\mu_{s+p+1}}\cdot {\dots} \cdot \sum\limits_{\mu_{\nu}=0}^{n_{\nu}-2} \mu_{\nu}\Delta^2 b_{\mu_{\nu}}.
\end {equation}

Using estimate \eqref{2.29} in \eqref{2.26} we get
\begin{equation}\label{2.30}
Q_{R}(n_1) = O\biggl( \prod\limits_{i=1}^{s}\sqrt{\log n_i}\,\cdot b_{n_i}\cdot
\prod\limits_{i=s+1}^{s+p} (n_i-1) \Delta^{1} b_{n_{i}-1} \cdot
\prod\limits_{i=s+p+1}^{\nu} \sum\limits_{\mu_{i}=0}^{n_{i}-2} \mu_{i}\Delta^2 b_{\mu_{i}}\biggl).
\end {equation}
Note that if $s+p>0$, then, in view of \eqref{2.5} and the definition of $b_i$ \eqref{2.4}, we obtain from \eqref{2.30}
\begin{equation}\label{2.31}
Q_{R}(n_1)= o(1) \quad \text{as} \quad n_1,\dots,n_{\nu} \to \infty,\quad \text{if } \quad 0\le p,s \le \nu,\, 0 < s+p \le \nu.
\end {equation}

Thus,
\begin{equation}\label{2.32}
Q_{R}(n_1)= o(1) \quad \text{as} \quad n_1,\dots,n_{\nu} \to \infty \quad \text{for any} \quad R\in \Omega.
\end {equation}
The similar estimate is true for $Q_{R}(\alpha_1)$ (considering the definition of $\alpha_1$):
\begin{equation}\label{2.33}
Q_{R}(\alpha_1)= o(1) \quad \text{as} \quad n_1,\dots,n_{\nu} \to \infty\quad \text{for any} \quad R\in \Omega.
\end {equation}

Consider now the case $s+p=0$, i.e., $s=p=0$. It means that $r_1=\dots=r_{\nu}=2$. By setting in \eqref{2.27} $s=p=0$, in account of notation \eqref{2.12}, we have
$$
Q_{(2,\dots,2)}(n_1)=\prod\limits_{i=1}^{\nu}\Delta^{2} b_{n_i-2}\, A_{n_1-2,\dots,n_{\nu}-2}^{(3,\dots,3)}=
$$
$$
=\sum\limits_{\mu_{1}=0}^{n_{1}-2} \Delta^2 b_{\mu_{1}}\sum\limits_{l_{1}=0}^{\mu_{1}}\ \dots \biggl\{\sum\limits_{\mu_{\nu}=0}^{n_{\nu}-2} \Delta^2 b_{\mu_{\nu}} \sum\limits_{l_{\nu}=0}^{\mu_{\nu}}
S_{l_1,\dots,l_{\nu},n'}(x;g)\biggl\}.
$$
In this case, by \eqref{2.19}, \eqref{2.23} and \eqref{2.25}, we obtain:
$$
Q_{(2,\dots,2)}(n_1)-Q_{(2,\dots,2)}(\alpha_1)=\Delta^{2} b_{n_1-2}\, \prod\limits_{i=2}^{\nu}\Delta^{2} b_{n_i-2}\, A_{n_1-2,n_2-2,\dots, n_{\nu}-2}^{(3,\dots,3)}
$$
$$
- \Delta^{2} b_{\alpha_1-2}\prod\limits_{i=2}^{\nu}\Delta^{2} b_{n_i-2}\, A_{\alpha_1-2,n_2-2,\dots,n_{\nu}-2}^{(3,\dots,3)}
$$
$$
=\sum\limits_{\mu_{1}=\alpha_1-1}^{n_{1}-2} \Delta^2 b_{\mu_{1}} \sum\limits_{\mu_{2}=0}^{n_{2}-2} \Delta^2 b_{\mu_{2}}\dots \sum\limits_{\mu_{\nu}=0}^{n_{\nu}-2} \Delta^2 b_{\mu_{\nu}}\biggl\{\sum\limits_{l_{1}=0}^{\mu_{1}}\dots\sum\limits_{l_{\nu}=0}^{\mu_{\nu}}
S_{l_1,\dots,l_{\nu},n'}(x;g)\biggl\}.
$$
By \eqref{2.5} and the definition of $\alpha_1$ we have:
$\sum\limits_{\mu_{1}=\alpha_{1}-1}^{n_{1}-2} \mu_{1}\Delta^2 b_{\mu_{1}}=o(1)$ as $n_{1}\to \infty $. Thus, by Corollary of Theorem 1 the same way as above we obtain that, as $n_{1},\dots, n_{\nu} \to \infty $
\begin{equation}\label{2.34}
Q_{(2,\dots,2)}(n_1)-Q_{(2,\dots,2)}(\alpha_1)=
O\biggl(\prod\limits_{i=2}^{\nu} \sum\limits_{\mu_{i}=0}^{n_{i}-2} \mu_{i}\Delta^2 b_{\mu_{i}}\cdot \sum\limits_{\mu_{1}=\alpha_{1}-1}^{n_{1}-2} \mu_{1}\Delta^2 b_{\mu_{1}}\biggl )=o(1).
\end {equation}

Thus, in view of \eqref{2.23},\eqref{2.25}, by estimates \eqref{2.32}-\eqref{2.34} we get that \eqref{2.8} is true for $j=1$.

Estimate \eqref{0.9} follows from \eqref{0.8} and results by E.Stein \cite{Stein} (moreover, even a slightly more strong estimate can be reduced from \eqref{0.8}, see, e.g. estimate (22) in \cite[p. 347]{Nikishin}).

Theorem 2 is proved.

\section{Proof of Theorem 3}\label{s2}

{\bf Proof of Theorem 3.}\,\,\footnote{\, In the proof of this theorem, some ideas represented in \cite{Sjolin} are used.} In order to prove the theorem, it is sufficient to prove the validity of estimate \eqref{0.12}; estimate \eqref{0.11} is deduced from it by means of standard argumentation.

Let us fix an arbitrary sample $J_{N-2} \subset M$. Without loss of generality, we consider that $J_{N-2}=\{1,\dots,N-2\}$, $M\setminus J_{N-2}=\{N-1,N\}$.

We introduce the following notations which permit to carry out the proof with less complexity. Let $t= m_{N-1}$, $q=m_N$ and
\begin{equation}\label{3.1}
n'=(n_1^{(\lambda_1)},\dots,n_{N-2}^{(\lambda_{N-2})}) \in\mathbb Z^{N-2}, \quad m'=(m_1,\dots,m_{N-2})\in\mathbb Z^{N-2}.
\end{equation}
Thus, the vectors $n^{(\lambda)}[J_{N-2}]\in {\mathbb Z}^N$ and $m\in {\mathbb Z}^N$ can be written in the form:
\begin{equation}\label{3.2}
n^{(\lambda)}[J_{N-2}]=(n',n_{N-1},n_N),\quad
m= (m', t, q).
\end{equation}
Denote also
\begin{equation}\label{3.3}
x'=(x_1,\dots, x_{N-2}) \in\mathbb T^{N-2}, \quad u'=(u_1,\dots, u_{N-2}) \in\mathbb T^{N-2}.
\end{equation}
We represent the partial sum $S_{n^{(\lambda)}[J_{N-2}]}(x;f)$ in the real form (considering notations \eqref{3.1}- \eqref{3.3}):
\begin{equation}\label{3.4}
S_{n^{(\lambda)}[J_{N-2}]}(x;f)=\sum\limits_{m'=0}^{n'}\biggl\{\sum\limits_{t=0}^{n_{N-1}} \sum\limits_{q=0}^{n_N}G_{(m',t,q)}(x,f)\biggl\}.
\end {equation}

Further denote
\begin{equation}\label{3.5}
l(t,\,q)=\bigr\{\log\bigl[\,\min(|t|,|q|)+2\bigr]\bigr\}^{-1}.
\end{equation}
Thus, the condition \eqref{0.10} in Theorem 3 looks as follows:
$$
\sum_{m\in {\mathbb Z}^N}|c_m(f)|^2 l^{-2}(t,\,q)<+\infty,
$$
and hence, according to the Riezc-Fischer theorem, there exists a function $g\in L_2(\mathbb T^N)$ such that the Fourier coefficients of the functions $f$ and $g$ are connected by relations
\begin{equation}\label{3.6}
c_m(f)=c_m(g)\,l(t,\,q),\quad m\in {\mathbb Z}^N.
\end {equation}
So, the partial sum in \eqref{3.4} can be rewritten as follows
\begin{equation}\label{3.7}
S_{n^{(\lambda)}[J_{N-2}]}(x;f)=\sum\limits_{m'=0}^{n'}\biggl\{\sum\limits_{t=0}^{n_{N-1}} \sum\limits_{q=0}^{n_N}G_{(m',t,q)}(x,g)\cdot l(t,q)\biggl\}.
\end {equation}

Denote as $G$ the sum in the braces in \eqref{3.7}, i.e.,
\begin{equation}\label{3.8}
G= \sum\limits_{t=0}^{n_{N-1}}\sum\limits_{q=0}^{n_N}G_{(m',t,q)}(x,g)\, l(t,q).
\end {equation}
Further, let us again introduce "shorthand"\, notations. Taking into account that $l(t,q)=l(q,t)$ (see \eqref{3.5}), we set
$$
\Delta_ql(t,q)=l(t,q)-l(t,q+1),\, \Delta_t(\Delta_ql(t,q))= \Delta_q l(t,q)- \Delta_q l(t+1,q)
$$
\begin{equation}\label{3.9}
= l(t,q)-l(t+1,q)-l(t,q+1)+l(t+1,q+1),
\end {equation}
and as well
\begin{equation}\label{3.10}
G_{t,q}= G_{(m',t,q)}(x,g); \quad U_i(\beta)=\sum\limits_{j=0}^{\beta}G_{i,j}.
\end {equation}
Applying the Abel transformation to the sum over $q$ in \eqref{3.8} and considering notations \eqref{3.9}, \eqref{3.10}, we have:
$$
G=\sum\limits_{t=0}^{n_{N-1}}\Biggl[\sum\limits_{q=0}^{n_N-1}\Delta_ql(t,q)U_t(q)+ l(t,n_N) U_t(n_N)\Biggr]
$$
\begin{equation}\label{3.11}
=\sum\limits_{q=0}^{n_N-1}\sum\limits_{t=0}^{n_{N-1}}\Delta_ql(t,q)U_t(q)
+ \sum\limits_{t=0}^{n_{N-1}}\,l(t,n_N) U_t(n_N).
\end {equation}
Applying the Abel transformation to each sum over $t$ in \eqref{3.11} and denoting (in account of \eqref{3.10})
\begin{equation}\label{3.12}
V(\alpha,\beta)= \sum\limits_{i=0}^{\alpha}\sum\limits_{j=0}^{\beta}G_{i,j}= \sum\limits_{i=0}^{\alpha}   U_i(\beta),
\end {equation}
we obtain
$$
G=\sum\limits_{q=0}^{n_{N}-1}\Biggl[\,\sum\limits_{t=0}^{n_{N-1}-1}
\Delta_t(\Delta_q l(t,q))V(t,q)\Biggr] + \sum\limits_{q=0}^{n_{N}-1} \Delta_q l(n_{N-1},q) V(n_{N-1},q)
$$
$$
+\sum\limits_{t=0}^{n_{N-1}-1} \Delta_t l(t, n_N) V(t,n_N) + l(n_{N-1},n_N) V(n_{N-1},n_N)=I^{(1)}_{m',\,n_{N-1},\,n_N}(x,g)
$$
\begin{equation}\label{3.13}
+I^{(2)}_{m',\,n_{N-1},\,n_N}(x,g)+
I^{(3)}_{m',\,n_{N-1},\,n_N}(x,g)+I^{(4)}_{m',\,n_{N-1},\,n_N}(x,g).
\end {equation}

Returning to \eqref{3.7} and taking account of \eqref{3.8}-\eqref{3.13}, we have:
\begin{equation}\label{3.14}
S_{n^{(\lambda)}[J_{N-2}]}(x;f)=\sum\limits_{m'=0}^{n'}\sum\limits_{j=1}^4 I^{(j)}_{m',\,n_{N-1},\,n_N}(x,g)
=\sum\limits_{j=1}^4 I^{(j)}_{n^{(\lambda)}[J_{N-2}]}(x,g).
\end {equation}

{\bf Lemma 2.} {\it The following estimates are true
\begin{equation}\label{3.15}
\biggl\|\sup\limits_{\lambda_1,\dots,\lambda_{N-2},n_{N-1},n_N>0}|I^{(j)}_{n^{(\lambda)}[J_{N-2}]}(x,g)|
\biggl\|_{L_2({\mathbb T}^N)}\leq C\|g\|_{L_2({\mathbb T}^N)}, \quad j=1,2,3,4.
\end {equation}}
{\bf Proof of Lemma 2.}
Note that in view of the definition of $l(t,q)$ -- \eqref{3.5} and the differences $\Delta_j$ -- \eqref{3.9}, we have
\begin{equation}\label{3.16}
\Delta_t(\Delta_ql(t,q))=0 \quad \text{if}\quad t\neq q.
\end {equation}
Denote $l(s)=l(s,\,s)$ and $\Delta l(s)=l(s)-l(s+1)$.

Let us prove estimate \eqref{3.15} for $j=1$. In account of \eqref{3.13}, as well as notations \eqref{3.10}, \eqref{3.12}, we have
$$
I^{(1)}_{n^{(\lambda)}[J_{N-2}]}(x,g)=\sum\limits_{m'=0}^{n'}I^{(1)}_{m',\,n_{N-1},\,n_N}(x,g)
=\sum\limits_{m'=0}^{n'}\sum\limits_{t=0}^{n_{N-1}-1}\sum\limits_{q=0}^{n_N-1}\Delta_t(\Delta_q l(t,q))V(t,q)
$$
$$
=\sum\limits_{m'=0}^{n'}\sum\limits_{t=0}^{n_{N-1}-1}\sum\limits_{q=0}^{n_N-1}\Delta_t(\Delta_q l(t,q)) \sum\limits_{i=0}^{t}\sum\limits_{j=0}^{q}G_{i,j}.
$$
From this, denoting as
\begin{equation}\label{3.17}
n_0=\min(n_{N-1},n_N),
\end {equation}
and considering \eqref{3.16} and \eqref{3.7}, we obtain
\begin{equation}\label{3.18}
I^{(1)}_{n^{(\lambda)}[J_{N-2}]}(x,g)=\sum\limits_{m'=0}^{n'}\sum\limits_{t=0}^{n_0-1}
\sum\limits_{i=0}^{t}\sum\limits_{j=0}^{t}G_{m',\,i,\,j}(x,g)\Delta l(t)=
\sum\limits_{t=0}^{n_0-1}S_{n',\,t,\,t}(x;g)\Delta l(t).
\end {equation}
Repeatedly applying the Abel transformation in \eqref{3.18} and taking into account argumentation in \cite[Ch. 13, Theorem (1.8)]{Zigmund2}, we obtain:
\begin{equation}\label{3.19}
\biggl|\sum\limits_{t=0}^{n_0-1}S_{n',\,t,\,t}(x;g)\Delta l(t)\biggr|\leq C\sup\limits_{t>0}\biggl|\frac {1}{t+1}\sum\limits_{i=0}^{t}S_{n',\,i,\,i}(x;g)\biggr|,\quad x\in{\mathbb T}^N.
\end {equation}

The following result is a particular case of the theorem proved by us (see \cite[Theorem 1]{Bloshanskii1}).

{\bf Theorem A.} {\it Let $k$, $2\le k\le N-1$, $N\geq3$, and the vector $(n_1^{(\lambda_1)},\dots,n_{N-k}^{(\lambda_{N-k})}$, $n_0, \dots,\, n_0) \in {\mathbb Z}^N_1$, where $\{n_i^{(\lambda_i)}\}$, $i=1,\dots,N-k$, are lacunary sequences, and $n_0\in {\mathbb Z}^1_1$. For any function $\varphi\in L_2({\mathbb T}^N)$ the estimate is true
\begin{equation}\label{3.20}
\biggl\|\sup\limits_{\lambda_1,\dots,\lambda_{N-k},n_0>0}
|S_{n_1^{(\lambda_1)},\dots,n_{N-k}^{(\lambda_{N-k})},\,n_0,\dots,\,n_0}(x;\varphi)|\biggl\|_{L_2({\mathbb T}^N)}\leq C\|\varphi\|_{L_2({\mathbb T}^N)},
\end {equation}
where the constant $C$ does not depend on the function $\varphi$.}

Applying in the right part of \eqref{3.19} estimate \eqref{3.20} with $k=2$, we get:
$$
\biggl\|\sup\limits_{\lambda_1,\dots,\lambda_{N-2},n_{N-1},n_N>0}\biggl|
\sum\limits_{t=0}^{n_0-1}S_{n',\,t,\,t}(x;g)\Delta l(t)\biggr|\biggl\|_{L_2({\mathbb T}^N)}\leq C\|g\|_{L_2({\mathbb T}^N)}.
$$
This estimate, in view of \eqref{3.18}, \eqref{3.19}, proves estimate \eqref{3.15} for $j=1$.

Let us prove estimate \eqref{3.15} for $j=2,3$. Consider $I^{(2)}_{n^{(\lambda)}[J_{N-2}]}(x,g)$. In account of \eqref{3.13}, \eqref{3.14} and notations \eqref{3.10}, \eqref{3.12}, we obtain:
$$
I^{(2)}_{n^{(\lambda)}[J_{N-2}]}(x,g)=\sum\limits_{m'=0}^{n'}I^{(2)}_{m',\,m_{N-1},\,m_N}(x,g) =\sum\limits_{m'=0}^{n'}\sum\limits_{q=0}^{n_N-1} \Delta_ql(n_{N-1},\,q) V(n_{N-1},q)
$$
\begin{equation}\label{3.21}
=\sum\limits_{m'=0}^{n'}\sum\limits_{q=0}^{n_N-1} \Delta_ql(n_{N-1},\,q)
\sum\limits_{t=0}^{n_{N-1}}U_t(q)
=\sum\limits_{m'=0}^{n'}\sum\limits_{t=0}^{n_{N-1}}\biggl\{\sum\limits_{q=0}^{n_N-1} \Delta_ql(n_{N-1},\,q)U_t(q)\biggl\}.
\end{equation}

Denote the expression in the braces as $\widetilde{I}$; given \eqref{3.9}, we have:
\begin{equation}\label{3.22}
\widetilde{I}=\sum\limits_{q=0}^{n_N-1} \Delta_ql(n_{N-1},\,q)U_t(q)=\sum\limits_{q=0}^{n_N-1}\{l(n_{N-1},\,q)-l(n_{N-1},\,q+1)\}U_t(q).
\end{equation}
Let us "simplify"\, $\widetilde{I}$; for this purpose consider two cases: $n_{N-1}>n_N$ and $n_{N-1}\leq n_N$. If $n_{N-1}>n_N$, then in the sum \eqref{3.22} $n_{N-1}>q$. Here, in account of the definition of $l(t,\,q)$ \eqref{3.5}, we get:
$$
\widetilde{I}=\sum\limits_{q=0}^{n_N-1}\{l(q)-l(q+1)\}U_t(q)=
\sum\limits_{q=0}^{n_0-1}\Delta l(q)U_t(q).
$$
Let now $n_{N-1}\leq n_N$, then
$$
\widetilde{I}=\sum\limits_{q=0}^{n_{N-1}}\{l(n_{N-1},\,q)-l(n_{N-1},\,q+1)\}U_t(q)+
\sum\limits_{q=n_{N-1}+1}^{n_N-1}\{l(n_{N-1},\,q)-l(n_{N-1},\,q+1)\}
$$
$$
\times U_t(q)=\sum\limits_{q=0}^{n_{N-1}-1}\{l(q)-l(q+1)\}U_t(q)=
\sum\limits_{q=0}^{n_0-1}\Delta l(q)U_t(q).
$$
In this case, from \eqref{3.21}, by \eqref{3.10} and \eqref{3.4}, we have:
$$
I^{(2)}_{n^{(\lambda)}[J_{N-2}]}(x,g)=
\sum\limits_{m'=0}^{n'}\sum\limits_{t=0}^{n_{N-1}}\sum\limits_{q=0}^{n_0-1}\Delta l(q)U_t(q)
$$
\begin{equation}\label{3.23}
=\sum\limits_{q=0}^{n_0-1}\sum\limits_{m'=0}^{n'}
\sum\limits_{t=0}^{n_{N-1}}\sum\limits_{j=0}^{q}G_{m',\,t,\,j}(x,g)\Delta l(q)=\sum\limits_{q=0}^{n_0-1}S_{n', n_{N-1},\,q}(x;g)\Delta l(q).
\end{equation}
The same way as above we can obtain
\begin{equation}\label{3.24}
I^{(3)}_{n^{(\lambda)}[J_{N-2}]}(x,g)=\sum\limits_{t=0}^{n_0-1}S_{n', t, n_N}(x;g)\Delta l(t).
\end{equation}
The same as for $j=1$ (see \eqref{3.19}), we again apply the Abel transformation in \eqref{3.23} and get: for $x\in {\mathbb T}^N$
$$
|I^{(2)}_{n^{(\lambda)}[J_{N-2}]}(x,g)|\leq C\sup\limits_{n_N>0}\biggl|\frac {1}{n_N+1}\sum\limits_{i=0}^{n_N}S_{n',\,n_{N-1},\,i}(x;g)\biggr|=
$$
considering the form of the Cezaro means
\begin{equation}\label{3.25}
=C\sup\limits_{n_N>0}\biggl|\int\limits_{{\mathbb T}^1}K_{n_N}(x_N-u_N) F(u_N,\circ)d u_N\biggr|=
\sigma_{n_N}(x_N;F),
\end {equation}
where $K_{n_N}(u)$ is the Fej\.{e}r kernel and
\begin{equation}\label{3.26}
F(u_N,\circ)= S_{n',\,n_{N-1}}(x',x_{N-1}; g; u_N).
\end {equation}

We apply in \eqref{3.25} estimate \eqref{1.13} and therefore get:
$$
\biggl\|\sup\limits_{\lambda_1,\dots,\lambda_{N-2}>\,0, \atop{n_{N-1},\,n_N>\,0}}|I^{(2)}_{n^{(\lambda)}[J_{N-2}]}(x,g)|
\biggl\|_{L_2({\mathbb T}^N)}\leq \biggl\|\sup\limits_{\lambda_1,\dots,\lambda_{N-2}>0,\atop{n_{N-1}>\,0}}\{\sup\limits_{n_N>0}|\sigma_{n_N}(x_N;F)|\}
\biggl\|_{L_2({\mathbb T}^N)}
$$
$$
\leq \biggl\|\sup\limits_{\lambda_1,\dots,\lambda_{N-2}>0, \atop{n_{N-1}>\,0}}|S_{n',\,n_{N-1}}(x',x_{N-1}; g; u_N)|
\biggl\|_{L_2({\mathbb T}^N)}.
$$
Applying inequality \eqref{0.6}, we obtain
$$
\biggl\|\sup\limits_{\lambda_1,\dots,\lambda_{N-2}>\,0, \atop{n_{N-1},\,n_N>\,0}}|I^{(2)}_{n^{(\lambda)}[J_{N-2}]}(x,g)|
\biggl\|_{L_2({\mathbb T}^N)}\leq C\|g\|_{L_2({\mathbb T}^N)}.
$$

Estimate \eqref{3.15} for $j=2$ is proved. Proof of this estimate for $j=3$ (see \eqref{3.24}) is similar.

And finally, let us prove estimate \eqref{3.15} for $j=4$. From \eqref{3.13}, taking into account \eqref{3.4}, \eqref{3.10}, \eqref{3.12}, we have:
$$
I^{(4)}_{n^{(\lambda)}[J_{N-2}]}(x,g)=\sum\limits_{m'=0}^{n'}I^{(4)}_{m',\,n_{N-1},\,n_N}(x,g)=
\sum\limits_{m'=0}^{n'}l(n_{N-1},\,n_N)\cdot V(n_{N-1},\,n_N)
$$
$$
=l(n_{N-1},\,n_N) \sum\limits_{m'=0}^{n'}
\sum\limits_{t=0}^{n_{N-1}}\sum\limits_{q=0}^{n_N}G_{m',\,t,\,q}(x,g)
$$
\begin{equation}\label{3.27}
=S_{n',\,n_{N-1},\,n_N}(x;g)\cdot l(n_{N-1},\,n_N)=S_{n^{(\lambda)}[J_{N-2}]}(x;g)\cdot l(n_{N-1},\,n_N).
\end {equation}

With the help of the function $F(u_{N},\circ)$, defined in
\eqref{3.26}, we represent the partial sum $S_{n^{(\lambda)}[J_{N-2}]}(x;g)$ in the form:
\begin{equation}\label{3.28}
S_{n^{(\lambda)}[J_{N-2}]}(x;g)=
\frac {1}{\pi}\int\limits_{{\mathbb T}^1}D_{n_{N}}(u_{N})F(x_{N}+u_{N},\circ)du_{N}.
\end {equation}
Further, using standart argumentation (see \cite[p. 84]{Sjolin}) and notation \eqref{3.3}, from \eqref{3.28} we obtain for $x=(x',x_{N-1},x_N)\in\mathbb T^N$:
$$
\sup\limits_{\lambda_j>\,0,\, j\in J_{N-2},\atop {n_{N-1}>\,0,\, n_N\geq2}}|S_{n^{(\lambda)}[J_{N-2}]}(x;g)|\{\log n_N\}^{-1}
$$
\begin{equation}\label{3.29}
\leq
C\cdot \mathfrak{M}\biggl\{\sup\limits_{\lambda_j>\,0,\, j\in J_{N-2},\atop {n_{N-1}>0}}
|S_{n',n_{N-1}}(x',x_{N-1},g;u_N)|\biggr\},
\end {equation}
where $\mathfrak{M}(\circ)$ is the Hardy-Littlewood maximal function. From \eqref{3.29}, using inequality \eqref{0.6}, we get:
$$
\biggl\|\sup\limits_{\lambda_j>0,\, j\in J_{N-2},\atop {n_{N-1}>0,\, n_N\geq 2}}|S_{n^{(\lambda)}[J_{N-2}]}(x;g)|\{\log n_N\}^{-1}\biggl\|_{L_2({\mathbb T}^N)} \leq
C\|g\|_{L_2({\mathbb T}^N)}.
$$
Analogously we can prove
$$
\biggl\|\sup\limits_{\lambda_j>0,\, j\in J_{N-2},\atop {n_{N-1}\geq2,\, n_N> 0}}|S_{n^{(\lambda)}[J_{N-2}]}(x;g)|\{\log n_{N-1}\}^{-1}\biggl\|_{L_2({\mathbb T}^N)} \leq
C\|g\|_{L_2({\mathbb T}^N)}.
$$
The last two inequalities and estimate \eqref{3.27} (in account of the definition of $l(n_{N-1},\,n_N)$  \eqref{3.5}) prove the validity of estimate \eqref{3.15} for $j=4$.

Lemma 2 is proved.

From estimates \eqref{3.14}, \eqref{3.15}, \eqref{3.5}, \eqref{3.6} the validity of estimate \eqref{0.10} follows.

Theorem 3 is proved.

 \vspace {0.5 cm}
\begin{centering}

\end{centering}

\begin{thebibliography}{90}

\bibitem{Bloshanskii}  I.L. Bloshanskii, Weyl multipliers and the growth of the partial sums of rectangulary summable multiple trigonometric Fourier series, Sov. Math. Dokl., 44 (3) (1992) 749-752.

\bibitem{Bloshanskii1}  I.L. Bloshanskii, D.A.Grafov, Sufficient conditions for convergence almost everywhere of multiple trigonometric Fourier series with lacunary sequence of partial sums, Real Analysis Exchange 41 (1) (2015/2016) 159-172.

\bibitem{Carleson} L. Carleson, On convergence and growth of partial sums of Fourier series, Acta Math. 116 (1966) 135–157.

\bibitem{Chen} Jau D. Chen, Narn Rueih Shieh,  On a Sufficient Condition for the Convergence of Multiple Fourier Series, Bull. Inst. Math. Acad. Sinica. 5 (2) (1977) 391-395.

\bibitem{Fefferman}  C. Fefferman, On the divergence of multiple Fourier series, Bull. Amer. Math. Soc. 77 (2) (1971) 191-195.

\bibitem{Hardy} G. Hardy, W.W. Rogosinski, Fourier series, Cambridge Univ. Press, 1946.

\bibitem{Kaczmarz} S. Kaczmarz, Zur Theorie der Fouriersche Doppelreihen, Stud. Math. 2 (1) (1930) 91-96.

\bibitem{Kojima}  M. Kojima, On the almost everywhere convergence of rectangular
partial sums of multiple Fourier series, Sci. Repts. Kanazawa
Univ. 22 (2) (1977) 163-177.

\bibitem{Kolmogorov2} A.N. Kolmogoroff, G.A.Seliverstoff, Sur la convergence des s\'{e}ries de Fourier, Atti Accad. naz. Lincei. Rend. 3 (1926) 307-310.

\bibitem{Moricz}  F. Moricz, Lebesgue functions and multiple function series. I, Acta Math. 37 (4) (1981) 481-496.

\bibitem{Nikishin} E.M. Nikishin, Weyl multipliers for multiple Fourier series, Math. USSR Sbornik 18 (1972) 351-360.

\bibitem{Plessner} A. Plessner, "\"{U}ber Konvergenz von trigonometrischen Reihen, Jornal f\"{u}r reine und angew. Math. 155 (1926) 15-25.

\bibitem{Sanadze} D. Sanadze, Sh. Kheladze, Convergence and divergence of multiple
Walsh- Fourier series, Tr. Tbilis. Mat. Inst. Razmadze 55 (1977) 93-106 [in Russian].

\bibitem{Sjolin} P. Sj\"{o}lin, Convergence almost everywhere of certain singular
integrals and multiple Fourier series, Arkiv Matem. 9 (1) (1971) 65-90.

\bibitem{Stein} E. Stein, On limits of sequences of operators, Ann. Math. 74 (1961) 140-170.

\bibitem{Tevzadze} N.R. Tevzadze, On the convergence of double Fourier series of square
summable functions, Soobsh. Acad. Nauk Gruzin. SSR. 58 (1970) 277-279.

\bibitem{Gigia} L.V. Zhizhiashvili, Some problems in the theory of simple and multiple trigonometric and ortogonal series, Russian Math. Surveys 28 (2) (1973).

\bibitem{Zigmund2} A. Zygmund, Trigonometric series, V.\:2. Cambridge Univ. Press, 1959.

\bibitem{Zigmund1} A. Zygmund, Trigonometric series, V.\:1. Cambridge Univ. Press, 1959.



\end{thebibliography}
\end{document}